\documentclass[pdflatex,final,10pt]{article}

\usepackage{framed,multirow}

\usepackage{pgf}
\usepackage{subfig}

\usepackage{amsmath,amstext,amssymb,amsthm}
\usepackage{bm}
\usepackage{psfrag}
\usepackage[off]{auto-pst-pdf}
\usepackage{graphicx}
\usepackage{empheq}
\usepackage{bbold}
\usepackage{color}
\usepackage{units}
\usepackage{showlabels}
\usepackage{url}
\usepackage{verbatim}
\usepackage[
          colorlinks=true,linkcolor=blue
          ,pdfproducer={Latex with hyperref}
          ,pdfcreator={pdflatex}
          ]{hyperref}

\title{Preserving the accuracy of numerical methods discretizing anisotropic elliptic problems}
\author{C. Yang$^{\dagger}$\footnote{Corresponding author}, F. Deluzet$^{\ddagger}$, J. Narski$^\ddagger$, \\[3em]
$^\ddagger$Universit\'e de Toulouse; UPS, INSA, UT1, UTM,\\ Institut de Math\'ematiques de Toulouse,\\
CNRS, Institut de Math\'ematiques de Toulouse UMR 5219,\\
F-31062 Toulouse, France,\\
name.forname@math.univ-toulouse.fr \\[1em]
$^\dagger$School of Mathematics,\\ Harbin Institute of Technology,\\ 92 West Dazhi Street, \\Nan Gang District, Harbin, 150001, China\\
yangchang@hit.edu.cn\\[1em]
$^\star$ Corresponding author.
}

\let\eps\varepsilon
\date{\today}

\newtheorem{remark}{Remark}[section]

\def\bb{\mathbf{b}}

\def\bn{\mathbf{n}}

\def\Id{\mathbb{I}\textrm{d}}
\def\VFT{\textrm{VF,2}}
\def\VFF{\textrm{VF,4}}
\def\VFP{\textrm{VF,p}}
\def\ASYM{\textrm{ASYM,2}}
\def\SYM{\textrm{SYM,2}}

\begin{document}

\maketitle

\begin{abstract}
  In this paper we study the loss of precision of numerical methods discretizing anisotropic problems and propose alternative approaches free from this drawback. The deterioration of the accuracy is observed when the coordinates and the mesh are unrelated to the anisotropy direction. While this issue is commonly addressed by increasing the scheme approximation order, we demonstrate that, though the gains are evident, the precision of these numerical methods remain far from optimal and limited to moderate anisotropy strengths. This is analysed and explained by an amplification of the approximation error related to the anisotropy strength. We propose an approach consisting in the introduction of an auxiliary variable aimed at removing the amplification of the discretization error. By this means the precision of the numerical approximation is demonstrated to be independent of the anisotropy strength.
\end{abstract}

\paragraph{Keywords}
Anisotropic equation, Plasma Physics, Asymptotic-Preserving schemes.

\maketitle

\section{Introduction}
\label{sec:introduction}

The physics of magnetized plasma is governed by anisotropic equations due to the large particle mobility along the magnetic field lines in comparison to the mobilities in the directions perpendicular to this field. Modelling these systems on larger scales gives rise to equations with a large diffusion coefficients along the anisotropy direction as stated by the model problem:
\begin{subequations}\label{sys:model}
\begin{gather}\label{eq:def:P}
  \left\{ 
    \begin{array}{ll}
      \displaystyle - \Delta_\perp \phi^\eps + \frac{1}{\eps} \Delta_\parallel \phi^\eps = f^\eps     & \text{ in } 
      \Omega, \\[3mm]
      \displaystyle \bn \cdot \left( \nabla_\perp \phi^\eps + \frac{1}{\eps} \nabla_\parallel \phi^\eps\right) = 0
      & \text{ on } \Gamma_N\,,\\[3mm]
      \phi^{\varepsilon }= 0
      & \text{ on } \Gamma_D\,,
    \end{array}
  \right.
\end{gather}
where $\bb$ denotes the vector field providing the direction of the magnetic field, $\bb$ verifying $\|\bb\|_2=1$, and the reciprocal of the asymptotic parameter $\eps$ defines the anisotropy strength. The parallel and perpendicular operators (with respect to the anisotropy direction $\bb$) are defined as
\begin{alignat}{5}
\nabla_\perp \psi &= \left(\Id - \bb \otimes \bb\right) \nabla \psi\,,& \qquad
\nabla_\parallel \psi &&= \left(\bb \otimes \bb\right) \nabla \psi\,,&  \\
\Delta_\perp \psi&= \nabla \cdot \left(\nabla_\perp \psi\right)\,,& \qquad
\Delta_\parallel \psi&&= \nabla \cdot \left(\nabla_\parallel \psi\right)\,,&
\end{alignat}
for any smooth function $\psi$, with $\Id$ the identity matrix and $\otimes$ the tensor product. The outward normal to the domain $\Omega$ is denoted $\bn$, $\Gamma_N\cup \Gamma_D$ are the domain boundaries, with $\bb\cdot \bn = 0$ on $\Gamma_D$ and $\bb\cdot \bn \neq 0$ on $\Gamma_N$. The flux associated to the model problem is denoted $\mathcal{Q}^\eps$ and defined by
\begin{equation}
\mathcal{Q}^\eps = \nabla_\perp \phi^\eps + \frac{1}{\eps} \nabla_\parallel \phi^\eps \,.
\end{equation}
\end{subequations}

Different difficulties are identified in the literature, regarding the numerical approximation of such problems. One of these difficulties is related to the deterioration of the condition number of matrices stemming from the discretization of these problems. This issue depends on the boundary conditions considered at each end of the magnetic field lines. In the model problem stated by Eq.~\eqref{eq:def:P} the parallel operator supplemented with the boundary condition on $\Gamma_N$ has a kernel containing all the functions with no gradients parallel to the $\bb$-field. The matrices issued from discretizations of this problem become consistent with a problem admitting an infinite amount of solutions as $\eps\to 0$, exhibiting therefore a condition number increasing with the anisotropy strength. We refer for instance to \cite{yang_iterative_2018} for an analysis of the condition number of these matrices.

This difficulty is mainly addressed by Asymptotic-Preserving methods \cite{degond_asymptotic_2009,degond_asymptotic-preserving_2012,narski_asymptotic_2013,chacon_asymptotic-preserving_2014,tang_asymptotic_2017,wang_uniformly_2018,deluzet_two_2019} restoring uniqueness in the limit $\eps\to 0$.

A second difficulty is also largely referred in the literature. It concerns the loss of accuracy of numerical approximations when the mesh is misaligned with the $\bb$-field. This issue is not necessarily related to the deterioration of the system matrix condition number. However, the numerical methods are observed to produce approximations with a poor precision, eventually meaningless for large anisotropy strengths. This is referred to as precision pollution (see for instance \cite{gunter_modelling_2005,van_es_finite-difference_2014,crouseilles_comparison_2015,ratnani_anisotropic_2016}) and can be explained by an amplification of the tuncation error of the parallel operator due to the heterogeneity of the diffusion coefficients. Different approaches are proposed in order to alleviate this pollution. The main idea is to increase the approximation order of the numerical methods in order to decrease the truncation error originating from the parallel operator discretization. This is for instance the path investigated in \cite{gunter_modelling_2005,gunter_finite_2007,van_es_finite-difference_2014,crouseilles_comparison_2015}. Field aligned reconstructions \cite{hariri_flux-coordinate_2013} make use of the weakness of the solution gradients along the magnetic field lines to interpolate a precise reconstruction of the parallel discrete derivatives. Similar ideas are also proposed in \cite{van_es_finite-difference_2014,van_es_finite-volume_2016} implementing a tracking of the field lines.

In the present paper, we propose a different approach based on a rescaling of the parallel gradients, the purpose being to vanish the gap between the parallel and perpendicular diffusion coefficients as presented in Eq.~\eqref{eq:def:P}. By this means, the cause of the pollution is removed rather than diminished, preventing any amplification of the truncation error and restoring an unaltered precision for the numerical method. This is achieved by the introduction of an auxiliary variable aiming at deriving an equation with both the parallel and the perpendicular operators at the same scale. This is an idea implemented in some Asymptotic-Preserving methods (see \cite{degond_asymptotic-preserving_2012,deluzet_two_2019}). The benefits of these techniques have already been emphasized regarding the condition number of the matrices issued from these approaches. In this paper, the purpose is to unravel new properties regarding the precision of the numerical approximation carried out thanks to a system with a rescaled parallel dynamic. In particular, we demonstrate that despite an accurate approximation of the solution, in the sense of the $H^1$-norm, a discretization of the flux $\mathcal{Q}^\eps$ cannot be reconstructed directly from the approximations of the solution derivatives. This feature is analysed and explained by the difficulty to provide an accurate approximation of the parallel gradients of the solution. A new reconstruction of the flux, implementing a rescaled parallel dynamic, is introduced in this paper. The discrete approximations of this quantity, free from any pollution of the precision is proposed, with an accuracy unrelated to the anisotropy strength, contrariwise to the flux 
classically reconstructed from the discrete derivatives of the solution approximation.

The outline of the paper is the following. The so-called precision pollution of anisotropic problems discretized on misaligned meshes is studied in Sec.~\ref{sec:pollutions}. These investigations are conducted for the finite difference and finite volume methods introduced in \cite{gunter_modelling_2005} and \cite{crouseilles_comparison_2015}. The discrepancy of the numerical method accuracy (as those used in \cite{crouseilles_comparison_2015}) is clearly emphasized thanks to the analysis of the associated truncation errors. To give a global picture of the numerical issues stemming from the discretization of anisotropic problems, we also characterize the impact of the limited computer arithmetic precision on the accuracy of numerical approximations. 
The principles of the parallel dynamic rescaling are introduced in Sec.~\ref{sec:rescaling}. The benefits on the precision pollution are then demonstrated. Together with this property, the need to reconstruct the parallel gradients from the auxiliary variable are motivated. A new reconstruction of the flux implementing a rescaled parallel dynamic is therefore introduced.
The advantages of numerical methods based on a rescaling of the parallel dynamic are numerically investigated in Sec.~\ref{sec:simu} and compared to discretizations usually harnessed for this class of problems with different frameworks: finite difference, finite volumes as well as finite elements and various approximation order (from second to seventh). The precision of the methods proposed herein are shown to be unaltered by the anisotropy strength, contrariwise to discretizations of the anisotropic problem despite the use of high order approximations.

\section{On the pollution of the scheme precision} \label{sec:pollutions}
\subsection{Solution manufacturing}\label{sec:solution:manufacturing}
In this section, the loss of precision of the numerical methods is related to the anisotropy of the coefficients multiplying the parallel or the perpendicular operators in the problem \eqref{eq:def:P}. 
To provide a quantitative analysis of the interplay between the parallel and perpendicular dynamics, the process of the solution manufacturing is implemented in the simplified context of an oblique homogeneous anisotropy direction. This framework is also considered to carry out the truncation error of the discretizations in Appendix~\ref{Appendix:Tronc}. The parallel dynamic refers to the component of the solution gradients aligned with the $\bb$-field. The perpendicular dynamic is the complementary component of the gradient, perpendicular to $\bb$.

We consider $\bb=(\alpha,\beta)^T$, with $\alpha^2+\beta^2=1$, together with the adapted coordinates $(X,Z)$ defined as
\begin{alignat}{5}\label{eq:def:aligned:Coordinates}
 X &= \alpha x + \beta y \,,&\qquad 
 Y &&= - \beta x + \alpha y &\,,
\end{alignat}
$(x,y)$ being the Cartesian coordinates. The $X$ coordinate, as constructed by Eq.~\eqref{eq:def:aligned:Coordinates}, is aligned with the anisotropy direction. We also
recall the expression the flux associated to the elliptic equation in the problem~\eqref{eq:def:P} 
\begin{equation}\label{eq:def:Qeps}
\mathcal{Q}^\eps = \left(\begin{array}{c} \mathcal{Q}^\eps_x \\ \mathcal{Q}^\eps_y \end{array}\right) = \left( \begin{array}{c}
\displaystyle \left( 1 - \alpha^2(1-\frac{1}{\eps}) \right)\frac{\partial \phi^\eps}{\partial x} - \alpha \beta(1-\frac{1}{\eps}) \frac{\partial \phi^\eps}{\partial y} \\
\displaystyle  - \alpha \beta(1-\frac{1}{\eps})\frac{\partial \phi^\eps}{\partial x} + \left( 1 - \beta^2(1-\frac{1}{\eps}) \right)\frac{\partial \phi^\eps}{\partial y} \\
\end{array}\right)
\end{equation}

We now introduce the function
\begin{equation}
\phi^{\eta}(x,z) = \tilde\phi_\perp(Y) + \eta\tilde\phi_\parallel(X)=\phi_\perp(x,y)+ \eta\phi_\parallel(x,y) \,,
\end{equation}
parametrized by the constant $\eta$. The two functions $\phi_\perp$ and $\phi_\parallel$ are assumed to be of magnitude one:
\begin{equation}\label{eq:def:sol:manufacture}
| \phi_\perp(x,y)| \sim |\phi_\parallel(x,y)| \sim 1\,.
\end{equation}
 The component $\phi_\perp$ defines the variations of $\phi^\eta$ in the directions perpendicular to $\bb$, with $\nabla_\parallel \phi_\perp = 0$, while $\phi_\parallel$ relates the parallel dynamics ($\nabla_\perp \phi_\parallel = 0$). Therefore, the parameter $\eta$ may be interpreted as the magnitude of the parallel gradients with respect to that of the perpendicular ones.
\paragraph{Unscaled parallel dynamics:}
The parallel and perpendicular variations of the function $\phi^\eta$ may be assumed to be the same order of magnitude which amounts to setting $\eta=1$. Inserting this definition into the model problem we have the following definition of the source term
\begin{equation*}
f^\eps = -\Delta_\perp \phi_\perp - \frac{1}{\eps}\Delta_\parallel \phi_\parallel \,.
\end{equation*}
Therefore, $f^\eps$ is not bounded when $\eps\to 0$. This choice may not be the most relevant, since it corresponds to both a flux $\mathcal{Q}^\eps$ and a source term $f^\eps$ unbounded in the limit of infinite anisotropies. 

\paragraph{Physics compatible solutions with isotropic fluxes}
This class of solutions are those defining a source term $f^\eps$ bounded irrespective of $\eps$. This requires that the parallel gradients of the function are small compared to that of the perpendicular directions. This property is harnessed in the field aligned reconstructions used for instance in \cite{hariri_flux-coordinate_2013}.

This condition is met for $\eta\leq\eps$, both the source term $f^\eps$ and the flux $\mathcal{Q}^\eps$ remaining bounded in the limit $\eps\to 0$. In this regime, the fluctuation of the solution along the anisotropy direction are much smaller than the variations in the perpendicular direction. The intense diffusion along the vector $\bb$ prevents the development of parallel gradients on scale larger than $\eps$.  The component $\phi_\perp$ carries the macroscopic variations of the solution, while $\phi_\parallel$ induces microscopic corrections in $\phi^\eta$. A precise approximation of $\phi_\parallel$ is therefore not mandatory to define an accurate approximation of the solution: capturing $\phi_\perp$ may be sufficient to have a good approximation of $\phi^\eta$ in the $H^1$-norm. However, the component $\phi_\parallel$ is significant in the definition of the flux $\mathcal{Q}^\eps = \nabla_\perp \phi + (1/\eps)\nabla_\parallel \phi=\nabla_\perp \phi_\perp + (\eta/\eps)\nabla_\parallel \phi_\parallel$ for values $\eta \sim \eps$.
The regime of interest is therefore the one defined by
\begin{equation}
\eta \sim \eps \,, \qquad \eps \ll 1 \,.
\end{equation}{}

\paragraph{Physics compatible solutions with anisotropic fluxes}
An intermediate regime may be identified, with solutions defined as
\begin{equation}\label{eq:sol:sqrt:eps}
  \phi^{\eta}(x,y) = \tilde\phi_\perp(Y) + \tilde\phi_\parallel({\eta} X)=\phi_\perp(x,y)+ \phi_\parallel({\eta} x,{\eta} y) \,,
\end{equation}
and $\eta \sim \sqrt{\eps}$, $\eps\ll 1$.
The source term $f^\eps$ derived from this definition remains bounded in the limit $\eps\to 0$, however the associated flux $\mathcal{Q}^\eps = \nabla_\perp \phi_\perp + ({1}/{\sqrt{\eps}}) \nabla_\parallel \phi_\parallel$ may not be bounded irrespective of $\eps$-values. The flux $\mathcal{Q}^\eps$ is anisotropic with a parallel component much larger than the perpendicular one. 

Note that in inserting the scaling relation $\eta \sim \eps$ in Eq.~\eqref{eq:sol:sqrt:eps} provides a solution similar to the one derived from Eq.~\eqref{eq:def:sol:manufacture} both giving rise to an isotropic flux $\mathcal{Q}^\eps$.

\subsection{Amplification of the approximation error and pollution of the scheme precision}\label{sec:pollution:para}

\subsubsection{Pollution of the reconstructed parallel flux}\label{sec:pollution:para:flux}
The reconstruction of the parallel dynamics is one of the difficulties characterizing this class of problems. This issue is manifest when the flux $\mathcal{Q}^\eps$ is recomposed from the derivatives of the solution with respect to the Cartesian coordinates ($\partial_x\phi$, $\partial_y \phi$). To illustrate more specifically these features let us consider the asymmetric scheme (see Eqs.~\eqref{eqs:Asym}) examined in \cite{gunter_modelling_2005}. This second order finite difference scheme provides the following approximations at $(x_{i+1/2},y_j)$:
\begin{subequations}\label{eq:def:dxh}
\begin{align}
\begin{split}
(\partial_x^\ASYM \Phi^h)_{i+1/2,j} &= \left(\frac{\partial \phi}{\partial x} \right) (x_{i+1/2},y_j) \\
& \hspace*{5em} +\frac{h^2}{24}\left(\frac{\partial^3 \phi}{\partial x^3} \right) (x_{i+1/2},y_j)+\mathcal{O}(h^4) \,,\label{eq:def:partialx}
\end{split}\\
\begin{split}
(\partial_y^\ASYM \Phi^h)_{i+1/2,j} &= \left(\frac{\partial \phi}{\partial y} \right) (x_{i+1/2},y_j)+\frac{h^2}{8}\left(\frac{\partial^3 \phi}{\partial y^3} \right) (x_{i+1/2},y_j)\\
& \qquad \qquad \quad+ \frac{h^2}{6}\left(\frac{\partial^3 \phi}{\partial x^2 \partial y} \right) (x_{i+1/2},y_j) + \mathcal{O}(h^4) \,.
\end{split}
\end{align}
\end{subequations}
These discrete operators are used to build the parallel gradient $\nabla_\parallel \phi = (\alpha \bb\cdot \nabla\phi, \beta \bb\cdot \nabla\phi)^T$
with the following expression for the first component $(\alpha \bb\cdot \nabla \phi)$ 
\begin{align}\label{eq:def:parallel:gradient}
\begin{split}
&\Big(\alpha \bb\cdot \nabla^\ASYM \Phi^h\Big)_{i+1/2,j} = (\alpha \bb\cdot \nabla\phi) (x_{i+1/2},y_j)-  \\
& \hspace*{5em}\frac{h^2}{24} \Bigg( \alpha^2 \frac{\partial^3\phi}{\partial x^3}(x_{i+1/2},y_j)+ \alpha\beta \Big( 3 \frac{\partial^3\phi}{\partial x^2\partial y}(x_{i+1/2},y_j)+ \\
& \hspace*{15em} 4 \frac{\partial^3\phi}{\partial y^3}(x_{i+1/2},y_j) \Big)\Bigg) +
 \mathcal{O}(h^4) \,,
\end{split}
\end{align}

To assess quantitatively the quality of the parallel gradient reconstruction, let us consider a specific definition of the solution with
\begin{equation}\label{def:sol:manufactured}
\phi^\eta=\phi_\perp(Y)+\eta \phi_\parallel(X)=\cos(2\pi n Y)+\eta \cos(2\pi m X)
\end{equation}
where $n$ and $m$ are two parameters with integer values. This amounts to assume that the solution is smooth and to perform its decomposition into Fourier modes. Then the discretization method is analysed for specific modes. Inserting this definition into \eqref{eq:def:parallel:gradient} yields
\begin{subequations}\label{subeq:error:GradPara:MN}
\begin{align}\label{eq:error:GradPara:M}
\begin{split}
&\left(\alpha \bb \cdot \nabla^\ASYM \Phi^h_\perp\right)_{i+1/2,j}=\\ 
& \hspace*{5em}-\frac{h^2}{3} (n\pi)^3 \, 2 \alpha^2 \beta \left(2\alpha^2+\beta^2\right)\sin(2\pi n Y_{i+1/2,j}) + \mathcal{O}(h^4) \,, 
\end{split}\\
\begin{split}\label{eq:error:GradPara:N}
&\left(\alpha \bb \cdot \nabla^\ASYM \Phi^h_\parallel\right)_{i+1/2,j}= -2 (m \pi) \sin(2 \pi mX_{i+1/2,j}) \\
  & \hspace*{3em} - \frac{h^2}{3} (m\pi)^3 \, \alpha\left(\alpha^4+3\beta^2\alpha^2+4\beta^4\right)\sin(2\pi m X_{i+1/2,j}) + \mathcal{O}(h^4) \,.
\end{split}
\end{align}
\end{subequations}

This finally provides, for the solution defined by Eq.~\eqref{def:sol:manufactured}
\begin{equation}{}\label{eq:error:GradPara:Final}
\begin{split}
&\left(\alpha \bb \cdot \nabla^\ASYM \Phi^{\eta,h}\right)_{i+1/2,j} = {\eta} \Bigg( -2 (m \pi) \sin(2 \pi mX_{i+1/2,j})
\\
&\hspace*{7em}-  \frac{h^2}{3} (m\pi)^3 \, \alpha\left(\alpha^4+3\beta^2\alpha^2+4\beta^4\right)\sin(2\pi m X_{i+1/2,j})\Bigg)  \\
& \hspace*{7em}- \frac{h^2}{3}  (n\pi)^3 \, 2 \alpha^2 \beta \left(2\alpha^2+\beta^2\right)\sin(2\pi n Y_{i+1/2,j}) + \mathcal{O}(h^4) \,.
\end{split}
\end{equation}

Some conclusions may be drawn from the estimates stated by Eqs.~\eqref{eq:error:GradPara:Final}. First, remark that, for any smooth function $\phi$, the following identities hold true
\begin{subequations}
\begin{alignat}{3}
\nabla_\parallel \phi &= 
  \left(\begin{array}{c}  \alpha_\parallel \partial_x \phi + \beta_\parallel \partial_y \phi \\
\gamma_\parallel \partial_x \phi + \kappa_\parallel \partial_y \phi
\end{array}\right)\,, \qquad 
\nabla_\perp \phi &&= \left(\begin{array}{c}
\alpha_\perp \partial_x \phi + \beta_\perp \partial_y \phi \\
\gamma_\perp \partial_x \phi + \kappa_\perp \partial_y \phi
\end{array}\right)\,;
\end{alignat}
where $\alpha_{\parallel,\perp}$, $\beta_{\parallel,\perp}$, $\gamma_{\parallel,\perp}$ and $\kappa_{\parallel,\perp}$ are related to the $\bb$-field coordinates, with the following properties for non-aligned coordinates
\begin{equation}
  0<|\alpha_{\parallel,\perp}|< 1 \,, \quad 0<|\beta_{\parallel,\perp}|< 1 \,, \quad 
  0<|\gamma_{\parallel,\perp}|< 1 \,, \quad 0<|\kappa_{\parallel,\perp}|< 1 \,.
\end{equation}
\end{subequations}

These relations mean that the condition
$\|\nabla_\perp \phi^\eta(x,y)\|_2 \gg \|\nabla_\parallel \phi^\eta(x,y)\|_2 $
is met if $\|\nabla \phi^\eta(x,y)\|_2\sim\|\nabla_\perp \phi^\eta(x,y)\|_2 \gg \|\nabla_\parallel \phi^\eta(x,y)\|_2$.
From these assertions, we may infer that, for scaled parallel dynamics ($\eta \leq \eps$) and large anisotropies ($\eta=\eps \ll 1$) the magnitude of the solution derivatives with respect to $x$ and $y$ are comparable to that of the components of the perpendicular gradients: $|\partial_{x,y} \phi^\eta(x,y)|\sim \|\nabla_\perp \phi^\eta(x,y)\|_2$. Therefore, the discrete parallel gradient cannot be reconstructed accurately using approximations of $\partial_{x,y} \phi^\eta$ for $\eps\ll1$.
This originates from the truncation error of the discretizations used for the solution derivatives as defined by Eqs.~\eqref{eq:def:dxh} which provide the estimates stated in Eq.~\eqref{eq:error:GradPara:M}. The truncation error for both the parallel gradient and the solution derivatives is proportional to $h^2$, while, considering a scaled parallel dynamic, the parallel gradient scale as $\eps$ and we recall $|\partial_{x,y} \phi^\eta|\sim\|\nabla_\perp \phi^\eta(x,y)\|_2 \sim 1$. The parallel gradients being rescaled by a factor $1/\eps$, the truncation error is amplified by this same ratio. From Eq.~\eqref{eq:error:GradPara:Final} the following scaling relation is stated for the component of the flux $\mathcal{Q}^\eps$ related to the parallel gradient:
\begin{equation}\label{eq:def:parallel:flux}
\begin{split}
&\left(\frac{\alpha}{\eps} \bb \cdot \nabla^\ASYM \Phi^{\eta,h}\right)_{i+1/2,j} \sim -2 (m \pi) \sin(2 \pi mX_{i+1/2,j})\\
& \hspace*{8em} - \frac{h^2}{\eps} \left( \frac{2}{3}  (n\pi)^3 \, \alpha^2 \beta \left(2\alpha^2+\beta^2\right)\sin(2\pi n Y_{i+1/2,j}) \right)\,.
\end{split}
\end{equation}
The meaningful contribution in this equation, proportional to $(m\pi)$, and is due to the parallel gradients of $\phi_\parallel$ while the error stemming from the perpendicular component (parallel gradients of $\phi_\perp$) has a magnitude proportional to $(n\pi)^3 h^2/\eps$. Therefore, the pollution stemming from the perpendicular dynamics discretization may be the dominant contribution in this equation with a magnitude scaling as $h^2/\eps$. This term deteriorates the precision of the numerical approximation for intermediate anisotropy strengths ($1>\eps > h^2$). For steepest anisotropies with $\eps < h^2$ the numerical method is ineffective.

Consider now the manufactured solution 
\begin{equation}{}\label{def:sol:manufactured:bis}
\phi^\eta=\phi_\perp(Y)+\phi_\parallel(\eta X)=\cos(2\pi n Y)+ \cos\big(2\pi m (\eta X) \big)
\end{equation}
with the following contribution to the reconstructed flux
\begin{equation}{}\label{eq:error:GradPara:Final:bis}
\begin{split}
& \left(\frac{\alpha}{\eps} \bb \cdot \nabla^\ASYM \Phi^{\eta,h}\right)_{i+1/2,j} \sim - \frac{\eta}{\eps}   \Big( 2 (m \pi) \sin(2 \pi m \eta X_{i+1/2,j})\Big)
\\
& \hspace*{8em}- \frac{h^2}{\eps} \left( \frac{2}{3}(n\pi)^3 \, \alpha^2 \beta \left(2\alpha^2+\beta^2\right)\sin(2\pi n Y_{i+1/2,j})\right) 
\end{split}
\end{equation}
We now investigate an intermediate regime, consisting of anisotropic fluxes but bounded source term. This amounts to set
$\eta \sim \sqrt{\eps}$. The pollution in this context is alleviated, the meaningful contribution in Eq.~\eqref{eq:error:GradPara:Final:bis} being offset when $h^2\sim \sqrt{\eps}$. This is to be compared to $h^2 \sim \eps$ for the isotropic flux case related by Eq.~\eqref{eq:def:parallel:flux}. It should be alose noted that, for both definitions, the deterioration of the precision increases with the magnitude of the solution perpendicular gradients (proportional to $n\pi$ for the examples defined by Eqs.~\eqref{def:sol:manufactured} and \eqref{def:sol:manufactured:bis}).

For unscaled parallel dynamics, which may be associated with the scaling relation $\eta =1$, the parallel and perpendicular gradients as well as the derivatives with respect to the Cartesian coordinates are comparable in magnitudes. Therefore, the pollution of the parallel gradient approximation does not occur in this context, with in the end, a reconstructed parallel dynamics at the right scale ($1/\eps$ and therefore not bounded when $\eps\to0$).

\subsubsection{Pollution of the discrete anisotropic equation}\label{sec:pollution:perp}
We now investigate a similar issue however examined under a different view point. The focus is now on the computation of the solution $\phi^\eps$ by means of a discretization of the anisotropic equation. To outline this specific feature, the solution is supposed to have vanishing parallel gradients. Therefore, inserting $\phi^\eps = \phi_\perp$ into the anisotropic equation, the only remaining contribution is $-\Delta_\perp \phi^\eps = -\Delta \phi_\perp$ owing to the identity $\Delta_\perp \psi= \Delta \psi - \Delta_\parallel\psi$.
This property holds true on the continuous level, however it is not exactly verified for the discrete quantities.   More specifically, assuming the following form for the macroscopic solution component
\begin{equation}\label{eq:def:phi:Perp}
\phi_\perp(x,y) = \cos\Big(2\pi n (-\beta x + \alpha y)\Big) \,,
\end{equation}
where $n\in \mathbb{N}$, yields the following truncation errors (see Appendix~\ref{Appendix:Tronc})
\begin{align}
(\Delta^\VFT_\parallel \bar\phi_\perp^h)_{i,j} &= -h^2 4 \left(n\pi\right)^4\,{(\alpha\, \beta)}^{2} \phi_\perp(x_i,y_j)  
+ \mathcal{O}(h^4) \,,\\
(\Delta^\VFF_\parallel \bar\phi_\perp^h)_{i,j} &= -h^4 \frac {323}{90} \left(n{\pi}\right)^{6} \alpha^2\beta^2\left( {\alpha}^{4}+{
\beta}^{4} \right) \phi_\perp(x_i,y_j) + \mathcal{O}(h^6) \,.
\end{align}
These quantities should be compared to 
\begin{align}
(\Delta^h_\perp \bar\phi_\perp^h)_{i,k} &= - 4 \left(n\pi\right)^2 \phi_\perp(x_i,y_j) + \mathcal{O}(h^p)\,,
\end{align}
to assemble the anisotropic differential operator applied to $\phi_\perp$. We indeed obtain:
\begin{subequations}\label{eq:Schemes:Aniso}
\begin{equation}\label{eq:equiv:VFT}
\begin{split}
&-(\Delta^\VFT_\perp \bar\phi_\perp^h+ \frac{1}{\eps} \Delta^\VFT_\parallel \bar\phi_\perp^h)_{i,j} = \\
& \hspace*{4em}  4 \left(n\pi\right)^2  \left( 1 + \frac{h^2}{\eps}  \left( n\pi \,\alpha \beta\right)^{2}\right) \phi_\perp(x_i,y_j)  + \mathcal{O}(h^2)+ \mathcal{O}\left(\frac{h^4}{\eps}\right) \,;
\end{split}
\end{equation}
for the second order approximations, with a similar identity for the fourth order one:
\begin{equation}\label{eq:equiv:VFF}
\begin{split}
&-(\Delta^\VFF_\perp \bar\phi_\perp^h - \frac{1}{\eps}\Delta^\VFF_\parallel \bar\phi_\perp^h)_{i,j} = 4 \left(n\pi\right)^2 \Bigg( 1\\
&\hspace*{2em}+ \frac{h^4}{\eps} \frac {323}{360} \left(n{\pi}\right)^{4} \alpha^2\beta^2\left( {\alpha}^{4}+{
\beta}^{4} \right)  \Bigg) \phi_\perp(x_i,y_j) + \mathcal{O}(h^4) + \mathcal{O}\left(\frac{h^6}{\eps}\right) \,.
\end{split}
\end{equation}
\end{subequations}
The mesh size is set to capture the derivatives of the function which, for scaled parallel dynamics, are comparable to the perpendicular gradients. Assuming that the mesh is refined to resolve precisely a period of $\phi_\perp$ (see Eq.~\eqref{eq:def:phi:Perp}) with 30 grid nodes, the following scaling relation can be stated
\begin{equation}
h \sim \frac{1}{30}\frac{1}{\pi n} \,,
\end{equation}

From \eqref{eq:equiv:VFT} it appears that the 
contribution of the perpendicular operator is totally offset by the discretization error of the parallel operator when $\eps<10^{-3}$ (for the second order discretizations). This threshold may be improved using higher order discretizations. Indeed with the fourth order scheme (see Eq.\eqref{eq:equiv:VFF}), the amplified truncation error is dominant when $\eps<10^{-6}$. Note also that the use of oversampled meshes, {\it i.e.} with increased values of $h \cdot (\pi n)$, offsets the influence of the amplified parallel discretization errors: the pollution is reduced when refining the mesh. However, the precision of the numerical method remains deteriorated. Though the convergence rate may be at the right order ($h^p$ for a discretization of order $p$) the precision is not optimal whatever the value of $\eps<1$: the precision of the numerical methods discretizing anisotropic problems with homogeneous $\bb$-fields is deteriorated compared to those discretizing isotropic problems.

\begin{remark}\label{remark:symetric:scheme}
The symmetric finite difference scheme defined by Eqs.~\eqref{eqs:Sym} \cite{gunter_modelling_2005} provides a truncation error 
\begin{equation}\label{eq:truncation:Sym:PhiPerp}
\begin{split}
(\Delta^\SYM_\parallel \phi_\perp^h )_{i,j} &= \frac{h^4}{2880} (n \pi)^6  \alpha^2 \beta^2 (\alpha-\beta)^2 (\alpha+\beta)^2\\
& \hspace*{8em}\big(h^2 (n\pi)^2-20 \big)  \phi_\perp(x_i,y_j) + \mathcal{O}(h^6)  \,;
\end{split}
\end{equation}
$(\phi_\perp^h)$ being the vector containing the values $\phi_\perp(x_i,y_j)$ (see Appendix~\ref{sec:appendix:notations}).
This discretization gives rise to a fourth order approximation of the parallel Laplacian applied to functions of the transverse coordinate ($Y$). This property holds true for homogeneous magnetic fields. However, for general (heterogeneous) anisotropies, the precision of this discretization remains second order accurate, comparable to the asymmetric scheme defined by Eqs.~\eqref{eqs:Asym} (See Appendix~\ref{sec:finite:difference}). This property will be highlighted in Sec.~\ref{sec:simu} dedicated to numerical investigations.
\end{remark}

To conclude this section, it is important to point out that for some specific frameworks the computation of an accurate approximation of the solution may be carried out. We can identify solutions with vanishing perpendicular gradients. Such solutions can be manufactured by considering the converse situation to the one analysed in the preceding lines, with a vanishing perpendicular component ($\phi_\perp$) yielding
\begin{equation}
\Delta_\perp \phi^\eta \ll \Delta_\parallel \phi^\eta 
\end{equation}
 For solutions only defined by their parallel component $\phi_\parallel$, the pollution does not alter the precision of the numerical approximation. A second class of the solution may also be considered: those with unscaled parallel dynamics as presented in Sec.~\ref{sec:solution:manufacturing}. Here also, the contribution of the perpendicular dynamics can be disregarded in the definition of the source term. Therefore polluting this contribution does not affect the quality of the numerical approximation.

\subsection{Impact of an arithmetic with a finite precision}\label{sec:arithmetic}
In this section, the deterioration of the precision of anisotropic problem discretizations related to the computer finite precision arithmetic is briefly addressed. In addition to the consistency issues outlined in Secs.~\ref{sec:pollution:para}, we emphasize that the round-off errors, due to a limited precision of the computer arithmetic, may also be amplified when discretizing anisotropic problems.

The impact of the arithmetic precision on the numerical approximation accuracy may be investigated in a simplified framework with an anisotropy aligned with one direction:
\begin{equation*}
\mathbf{b}=(1,0)^T.
\end{equation*}
Substituting this definition into~\eqref{eq:equiv:VFT} ou~\eqref{eq:equiv:VFF} with $\phi_\perp(y_j) = \cos(n\pi y_j)$ provides
\begin{align}\label{eq:SP:VFP}
-(\Delta^\VFP_\perp \bar\phi_\perp^h+ \frac{1}{\eps} \Delta^\VFP_\parallel \bar\phi_\perp^h)_{i,j} = 
4 \left(n\pi\right)^2  \phi_\perp(y_j)  + \mathcal{O}(h^p),\quad p\in\{2,4\}.
\end{align}
In this context and contrariwise to  $\bb$-fields misaligned with the mesh, the discrete parallel Laplacian vanishes when applied to any function of the transverse coordinate $y$. 
This means that the precision of the discrete perpendicular operator is not deteriorated by the discretization error of the parallel one. Therefore, a good precision should be obtained for any choice of $\eps$. 

To evaluate the influence of the inexact representation of numbers, let us introduce $\eps_A$ defined as
\begin{equation*}
\eps_A = \frac{|x-x_A|}{|x|}
\end{equation*}
the relative error of the number representation using the computer arithmetic, $x$ denoting the exact value and $x_A$ its representation in the computer arithmetic. For the double precision arithmetic the mean value of this parameter is usually equal to $10^{-16}$. The truncation error as stated by Eq.~\eqref{eq:SP:VFP} does not account for the round-off errors and should be corrected into
\begin{align}\label{eq:SP:VFP:roundoff}
-(\Delta^\VFP_\perp \bar\phi_\perp^h+ \frac{1}{\eps} \Delta^\VFP_\parallel \bar\phi_\perp^h)_{i,j} = 
4 \left(n\pi\right)^2  \phi_\perp(y_j)  + \mathcal{O}(h^p) +\mathcal{O}\left(\frac{\eps_A}{h^2\eps}\right) \,.
\end{align}
The additional error term on the right hand side of this equation represents the round-off errors due to the representation of the solution point values $\phi_{i,j}$ with the computer arithmetic, divided by the squared mesh size (due to the double discrete derivative) and amplified by the ratio $1/\eps$. A more subtle analysis should be carried to precisely evaluate the impact of the computer arithmetic on the precision of the problem solution. This would amount to substituting $1/(\eps h^2)$ by the condition number of the matrix associated with the discrete differential operators. However this heuristic provides a good insight of the issues related to the limited precision of the number representation.

The influence of the error made on the representation of floating point numbers may be disregarded when ${\eps_A}/({\eps h^2}) \ll h^p$, the accuracy of numerical approximations being driven by the precision of the discretizations. As the ratio ${\eps_A}/({\eps h^2})$ reaches the threshold $h^p$ the precision of the computations is altered and refining the mesh may deteriorate the accuracy rather than improving it. This is more penalizing for high order methods, the threshold $h^p$ being reached for coarser values of $h$.
Compared to isotropic problems, the propagation of round-off errors may be a serious issue since the amplification factor is proportional to the imbalance between the perpendicular and parallel dynamics ($1/\eps$).

\section{Preventing the loss of accuracy thanks to a rescaling of the parallel dynamic}\label{sec:rescaling}

\subsection{Introduction of a rescaled auxiliary variable}
The analyses carried out in the preceding section, unravel that the loss of accuracy is genuinely related to the stiffness of the equation due to the heterogeneity of the parallel and perpendicular diffusions. We propose here to raise this issue by rescaling the parallel dynamics thanks to the use of an auxiliary variable $q^\eps$ satisfying 
\begin{equation}
\nabla_\parallel \phi^\eps = \eps  \nabla_\parallel q^\eps \,,
\end{equation}
Injecting this identity into the anisotropic problem, the following system is derived
\begin{subequations}\label{subeq:MM}
\begin{align}
&- \Delta_\perp \phi^\eps - \Delta_\parallel q^\eps = f^\eps \,, \label{eq:MM:one}\\
&-\Delta_\parallel \phi^\eps = -\eps  \Delta_\parallel q^\eps \,. \label{eq:MM:two}
\end{align}
\end{subequations}
The heterogeneity of the parallel and perpendicular diffusions has disappeared in 
Eq.~\eqref{eq:MM:one} the multi-scale nature of the problem being contained in Eq.~\eqref{eq:MM:two} involving only the parallel dynamic.
Dicretizations of this problem are therefore not subjected to the pollution of the precision. This is actually a property of Asymptotic-Preserving schemes demonstrated within this work.
To shorten the presentation, only the discrete version of Eq.~\eqref{eq:MM:two} is analysed, since this is the only multi-scale equation. The finite volume discretization introduced in the section provides
\begin{subequations}\label{eq:discrete:MM:two}
\begin{equation}
\begin{aligned}\label{eq:discrete:MM:two:a}
&\left(\Delta_\parallel^\VFT (\bar \Phi^h- \eps \bar {\bm{q}}^h) \right)_{i,j} = \Delta_\parallel \phi(x_i,y_j) - \eps \Delta_\parallel q(x_i,y_j)\\
&\qquad+\alpha^2\frac{h^2}{12}\left( \frac{\partial^4}{\partial x^4} \Big((\phi - \eps q)(x_i,y_j) \Big) \right)\\
  &\qquad + \alpha\beta\frac{h^2}{3}\left(\frac{\partial^4}{\partial x^3 \partial y} \Big((\phi - \eps q)(x_i,y_j) \Big) + \frac{\partial^4}{\partial x \partial y^3} \Big((\phi - \eps q)(x_i,y_j) \Big)\right) \\
  &\qquad +\beta^2\frac{h^2}{12}\left( \frac{\partial^4}{\partial y^4} \Big((\phi - \eps q)(x_i,y_j) \Big)\right)+ \mathcal{O}(h^4)\,.
\end{aligned}
\end{equation}

\begin{equation}\label{eq:discrete:MM:two:b}
\begin{aligned}
&\left(\Delta_\parallel^\VFF (\bar \Phi^h- \eps \bar {\bm{q}}^h) \right)_{i,j} = \Delta_\parallel \phi(x_i,y_j) - \eps \Delta_\parallel q(x_i,y_j)\\
&\quad +\alpha^2 \frac{h^4}{1920}\left(\frac{\partial^6}{\partial x^2 \partial y^4} \Big((\phi - \eps q)(x_i,y_j) \Big) - \frac{64}{3}\frac{\partial^6}{\partial x^6} \Big((\phi - \eps q)(x_i,y_j) \Big) \right)\\
  &\quad - \alpha\beta\frac{h^4}{15}\left(\frac{\partial^6}{\partial x^5 \partial y} \Big((\phi - \eps q)(x_i,y_j) \Big) + \frac{\partial^6}{\partial x \partial y^5} \Big((\phi - \eps q)(x_i,y_j) \Big)\right) \\
  &\quad  +\beta^2 \frac{h^4}{1920}\Bigg(\frac{\partial^6}{\partial x^4 \partial y^2} \Big((\phi - \eps q)(x_i,y_j) \Big) \\
  & \hspace*{12em}- \frac{64}{3}\frac{\partial^6}{\partial y^6} \Big((\phi - \eps q)(x_i,y_j) \Big) \Bigg) + \mathcal{O}(h^6)\,.
\end{aligned}
\end{equation}
\end{subequations}
The discretizations of Eq.~\eqref{eq:MM:two} show a different picture compared to a straight discretization of the original anisotropic equation (see Eqs.~\eqref{eq:Schemes:Aniso}). When $\eps \ll h^p$ ($p$ depending of the precision order of the numerical method), the contribution of $q$ is lost in Eqs.~\eqref{eq:discrete:MM:two}. In the end we obtain a discretization of the equation
\begin{equation}\label{eq:MM:q:limit}
  - \Delta_\parallel \phi = 0 \,.
\end{equation}
This is actually a property satisfied by the solution of the problem when $\eps\to 0$ and the consistency with Eq.~\eqref{eq:MM:q:limit} translates the fact that Eqs.~\eqref{eq:discrete:MM:two} only account for the correction term  $\eps \Delta_\parallel q$ with the precision of the numerical scheme. 

Corollary, we cannot expect to recover a precise approximation of the parallel gradients ($\nabla_\parallel \phi$) from $\Phi^h$. Indeed, Eqs.~\eqref{eq:discrete:MM:two} cannot provide a reconstruction of $\nabla_\parallel \phi$ with a magnitude $\eps$ when $\eps \ll h^p$. This outlines the difficulty to account accurately of the parallel gradients for scaled parallel dynamics.  However, in this context, the variations of the problem solution are hardly explained by the parallel derivatives but the transverse ones ({\it i.e.} $\nabla_\perp \phi$). The perpendicular gradients appears in the Eq.~\eqref{eq:MM:one} where $\eps$ is expelled. Therefore, we can expect an accurate approximation of $\phi$ and $\nabla_\perp \phi$. In other words, a good approximation of $\phi$ in the $H^1$-norm can be anticipated. However, this is not sufficient to account for the parallel gradients accurately. Nonetheless, $\nabla_\parallel q$ is computed from Eq.~\eqref{eq:MM:one}
 with a precision comparable to that of $\nabla_\perp \phi$ . This means that both $\nabla_\perp \phi$ and $\nabla_\parallel q$ can be approximated with a precision independent of $\eps$. To take advantage of this property, we introduce the flux with a rescaled parallel dynamics $\mathcal{Q}^{RPD}$ defined as
\begin{equation}\label{eq:def:QRPD}
  \mathcal{Q}^{RPD} = \nabla_\perp \phi + \nabla_\parallel q \,.
\end{equation}
Both the parallel and the perpendicular gradients are inserted at the same scale in Eq.~\eqref{eq:def:QRPD}. The discretization of $\mathcal{Q}^{RPD}$ is therefore not subjected to the pollution of the scheme precision as stated in Sec.~\ref{sec:pollution:para:flux}.

\subsection{Implementations of the parallel dynamics rescaling}
The choice of the auxiliary variable $q$ introduced to rescale the parallel gradient is not unique. Two choices are presented here, the so-called Micro-Macro \cite{degond_asymptotic-preserving_2012} formulation and the Two-Field-Iterated method \cite{deluzet_two_2019,yang_numerical_2019}. 

The Micro-Macro method is the most straighforward implementation of the ideas introduced precedently to rescale the parallel dynamic thanks to an auxiliary variable. It corresponds to an auxiliary variable with a vanishing trace on one feet of each $\bb$-field line. 
The Micro-Macro formulation consists of the two coupled sets of eqations
\begin{subequations}\label{eq:def:MM} 
\begin{eqnarray}
&&\left\{
\begin{array}{lll}
-\Delta_{\perp} \phi - \nabla_\parallel q = f \,, & \text{in} & \Omega, \\ [1mm]
\mathbf{n}\cdot \left( \nabla_\perp \phi + \nabla_\parallel q \right) = 0, & \text{on} & \Gamma_N , \\[1mm]
\phi = 0,  & \text{on} & \Gamma_D , 
\end{array}
\right.\label{eq:def:MM:a} \\ [1mm]
&&\left\{
\begin{array}{lll}
-\Delta_\parallel \phi =  -\eps \Delta_\parallel q, & \text{in} & \Omega, \\ [1mm]
\mathbf{n}\cdot \nabla_\parallel \phi = \eps \mathbf{n}\cdot \nabla_{\parallel}q \,, & \text{on} & \Gamma_N , \\ [1mm]
q = 0,  & \text{on} & \Gamma_D \cup \Gamma_\text{in} \,,
\end{array}
\right.\label{eq:def:MM:b}
\end{eqnarray}
\end{subequations}
where
\begin{equation*}
\Gamma_\text{in} = \left\{ \mathbf{x}\in \Gamma_N \,|  \, \bb(\mathbf{x})\cdot \mathbf{n}(\mathbf{x})>0 \right\} \,.
\end{equation*}
The "inflow" condition prescribing a zero trace for $q$ on $\Gamma_\text{in}$ is mandatory to provide uniqueness of the auxiliary variable in the system \eqref{eq:def:MM}. However this condition requires that all the $\bb$-field lines intersect one boundary of the domain. Therefore, the Micro-Macro formulation is not well suited to address topologies including closed field lines.

To overcome this difficult, a different implementation of the rescaling may be proposed. The Two-Field Iterated method consists of a sequence $(\phi^n,q^n)$ satisfying the following equations
\begin{subequations}\label{eq:def:TFI}
\begin{eqnarray}
&&\left\{
\begin{array}{lll}
-\Delta_{\varepsilon_0}\phi^{n+1} = \varepsilon_0 f + (\varepsilon_0 - \varepsilon) \Delta_{\|} q^n, & \text{in} & \Omega, \\ [3mm]
\mathbf{n}\cdot \mathbb{A}_{\varepsilon_0}\nabla \phi^{n+1} = -(\varepsilon_0 - \varepsilon) \mathbf{n}\cdot \nabla_{\|} q^n, & \text{on} & \Gamma_N , \\ [3mm]
\phi^{n+1} = 0,  & \text{on} & \Gamma_D , 
\end{array}
\right.\label{eq:TFI_phi} \\ [3mm]
&&\left\{
\begin{array}{lll}
-\Delta_{\varepsilon_0} q^{n+1} =  f + \Delta_\bot(\phi^{n+1} -  \varepsilon_0 q^n), & \text{in} & \Omega, \\ [3mm]
\mathbf{n}\cdot \mathbb{A}_{\varepsilon_0}\nabla q^{n+1} = - \mathbf{n}\cdot \nabla_{\bot}(\phi^{n+1} - \varepsilon_0 q^n), & \text{on} & \Gamma_N , \\ [3mm]
q^{n+1} = 0,  & \text{on} & \Gamma_D .
\end{array}
\right.\label{eq:TFI_q}
\end{eqnarray}
where 
\begin{equation}
\Delta_{\varepsilon_0} \phi = \eps_0 \Delta_\perp \phi + \Delta_\parallel \phi \,,\qquad \mathbb{A}_{\varepsilon_0} = \eps_0 \left(\Id - \bb \otimes \bb \right) + \left(\bb \otimes \bb \right) \,
\end{equation}
\end{subequations}
In these equations, $\eps_0$ is a numerical parameter satisfying $1 > \eps_0 \gg \eps$.
The sequence is proved to converge towards $(\phi, q^\infty)$ (see \cite{deluzet_two_2019}), with $(\phi, q^\infty)$ satisfying the set of Eqs.~\eqref{eq:def:MM} but the inflow condition. Actually $q^\infty$ and $q$ differ by a function with no aligned gradients (in the kernel of the parallel operator). It is important to note that, though the sequence $(\phi^k,q^k)$ is constructed thanks to the resolution of mildly anisotropic problems, the fixed point $(\phi, q^\infty)$ is the solution of a system with no stiffness, therefore free from the pollution of the scheme precision.

\section{Numerical assessment of discretizations implementing a rescaled parallel dynamic}\label{sec:simu}

\subsection{Set-up definition}

The numerical method effectiveness is investigated by comparing the numerical approximations against manufactured solutions. These solutions are analytically derived from an exact (analytic) expression of the the $\bb$-field and the solution $\phi^\eps$. These expressions are used to compute analytically the source term of the anisotropic equation:
\begin{equation*}
  f^\varepsilon  = - \Delta_\perp \phi^{\varepsilon} - \frac{1}{\varepsilon }\Delta_\parallel \phi^{\varepsilon } \,.
\end{equation*}
The computational domain $\Omega$ is restricted to $[0,1]^2$.

The $\bb$-field components are defined as functions of the coordinates with
\begin{subequations}\label{eq:def:setup}
\begin{gather}
  \bb(x,y) = \frac{B(x,y)}{\|B(x,y)\|_2}\, , \quad
  B(x,y) =
  \left(
    \begin{array}{c}
      \theta  (2x-1) \cos (m\pi y) + \pi \\
      \pi \theta  m(x^2-x) \sin (m\pi y)
    \end{array}
  \right)
  \label{eq:J99a}\,,
\end{gather}
where $m$ and $\theta $ parametrize the topology of the field. For $\theta
=0$ the $\bb$-field reduces to $(1,0)^T$. It is aligned with the $x$-coordinate. For $0<\theta<\pi $
the field oscillates in the domain with $m/2$ periods and contains
open field lines only. For $\theta >\pi$ the field oscillates in the
domain and contains $m$ regions composed of closed field lines,
related to as magnetic islands in the context of plasma physics \cite{narski_asymptotic_2013} ---
see Fig.\ref{fig:test1_ana} for three configurations.
\def\xxx{0.3\textwidth}
\begin{figure}[!ht] 
  \centering{}
  \begin{tabular}{ccc}
    \resizebox{\xxx}{!}{\rotatebox{0}{\includegraphics{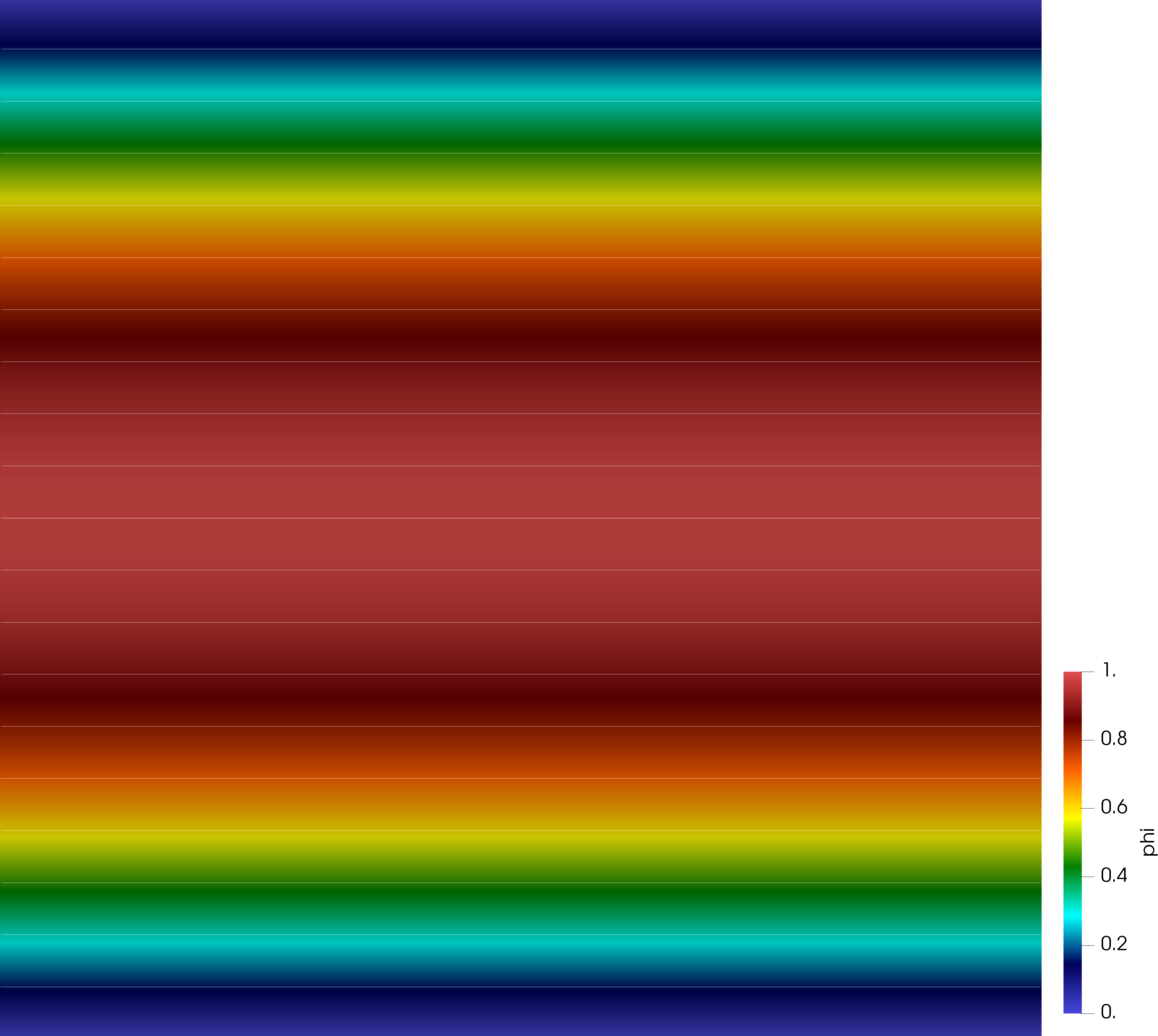}}}
    &
    \resizebox{\xxx}{!}{\rotatebox{0}{\includegraphics{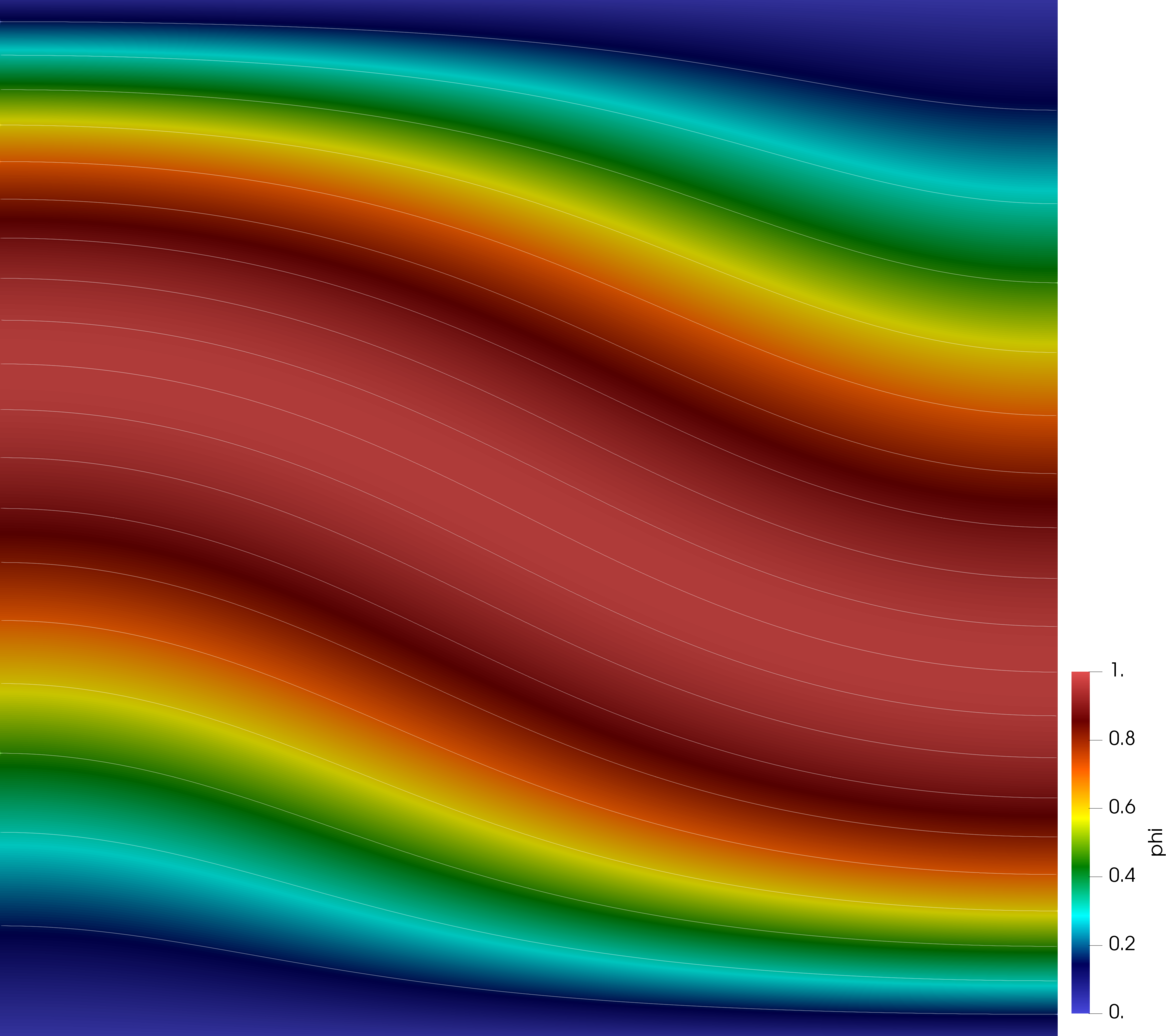}}}
    &
    \resizebox{\xxx}{!}{\rotatebox{0}{\includegraphics{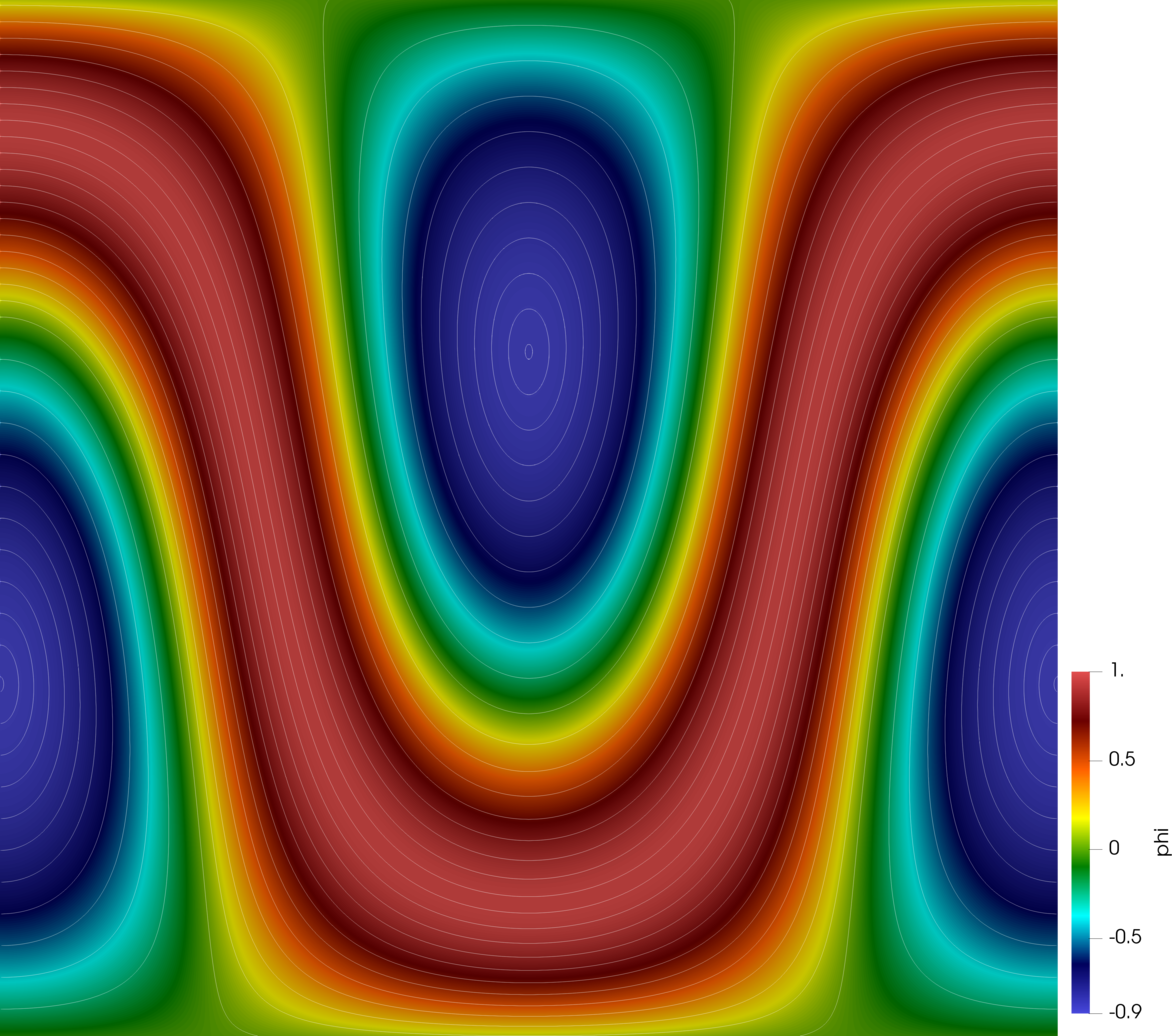}}}
    \\
    $\theta =0$
    &
    $\theta =2$, $m=1$
    &
    $\theta =10$, $m=2$
  \end{tabular}
\caption{Exact (manufactured) solution defined by Eq.~\eqref{eq:J79a} of the anisotropic problem \eqref{eq:def:P} in a colour scale for three sets of parameters defining the anisotropy direction specified by Eqs.~\eqref{eq:J99a}.}
\label{fig:test1_ana}
\end{figure}

Numerical tests are performed on manufactured solutions obtained by
adding a perturbation proportional to $\varepsilon $ to a function
constant in the direction of $\bb$:
\begin{gather}
  \phi^{\varepsilon } = \sin \Big( \omega \big(\pi x +\theta (x^2-x)\cos (m\pi y)\big)
  \Big)
  + \varepsilon
  \cos \left( 2\pi x\right) \sin \left(\pi    y \right) 
  \label{eq:J79a}\,.
\end{gather}
\end{subequations}
The parameter $\omega\in \mathbb{N}$ provides a control of the number of oscillations in the directions perpendicular to $\bb$. Its value is set to 1 except where specified.

Note that with the choice of $\phi ^\varepsilon $ stated by Eqs.~\eqref{eq:def:setup}, the framework is that of a scaled parallel dynamic: the source term is bounded independently of $\varepsilon $. 

\subsection{Numerical investigations}

\subsubsection{On the locking effect and the precision pollution}
In this section first numerical insights of the precision pollution are provided. This feature may be related to the so-called locking effect emphasized in the framework of Finite Element methods \cite{babuska_locking_1992}. The finite difference discretizations introduced in \cite{gunter_modelling_2005} are considered in this section with both the symmetric and the asymmetric schemes (see Sec.~\ref{sec:finite:difference}).

The locking effect is investigated on Fig.~\ref{fig:Locking:a} displaying the $L^2$-norm of the numerical approximation carried out by a finite difference discretization of the anisotropic problem. The computations are related to the set-up defined by Eqs.~\eqref{eq:def:setup}. With the increase of the anisotropy, the numerical approximation is observed to vanish (the $L^2$-norm of the solution decreases to 0), which is a typical feature of the locking effect \cite{babuska_locking_1992} also outlined in \cite{deluzet_two_2019}. Simultaneously, the precision of the numerical method is deteriorated. The computations carried out with a rescaled parallel dynamic are not subject to this loss of accuracy. 
\def\xxx{0.9\textwidth}
\begin{figure}[!ht]
\centering
\subfloat[\label{fig:Locking:a}Discrete anisotropic equation versus systems with a rescaled parallel dynamic.]{\resizebox{\xxx}{!}{\input{fig/Beta-Omega4-DF.tex}}{}}\newline\noindent
\subfloat[\label{fig:Locking:b}Corrected anisotropic equation as defined by Eqs.~\eqref{sys:def:aniso:corrected}.]{\resizebox{\xxx}{!}{\includegraphics[width=0.5\textwidth]{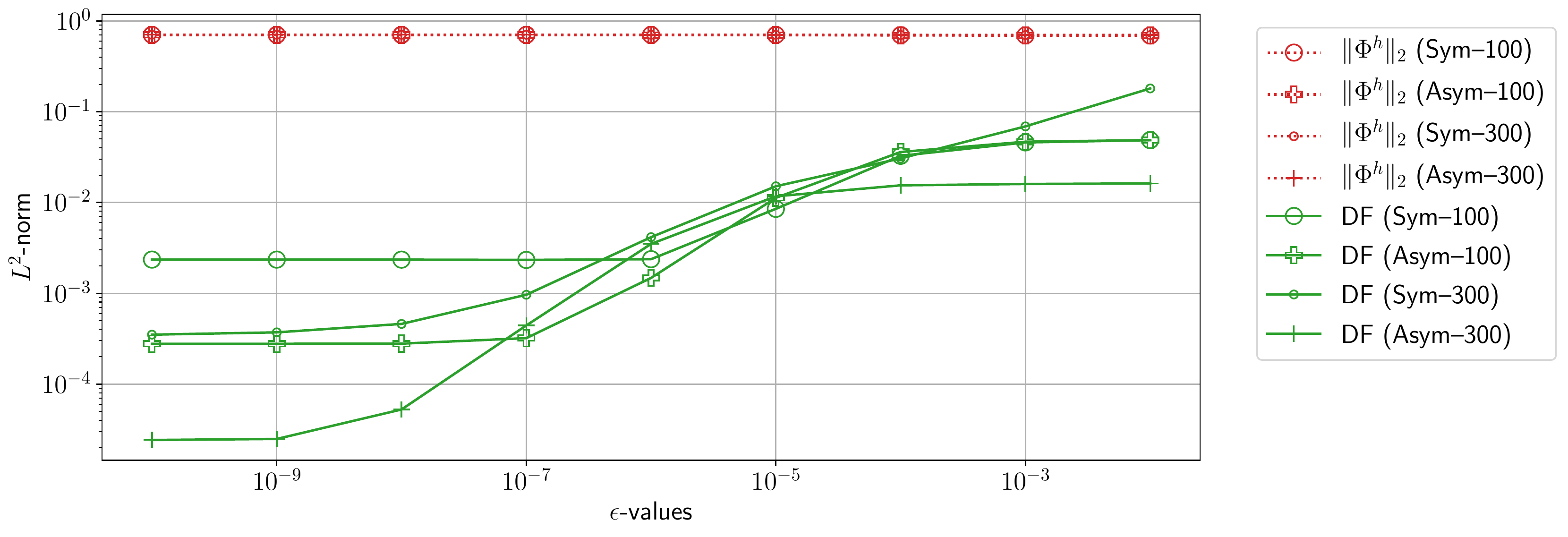}}{}}  
  \caption{Locking and pollution: $L^2$-norms of the numerical approximation $\|\Phi^h\|_2$ and of the absolute error against the exact solution for computations carried by a finite difference discretizations of the anisotropic problem (DF) and the system~\eqref{eq:def:TFI} with a rescaled parallel dynamics (RPD). The computations are related to the set-up~\eqref{eq:def:setup} with $\theta=2$, $m=1$, $\omega=4$ and carried out using the symmetric (Sym) or the asymmetric (Asym) scheme, on a $100\times100$ and a $300\times300$ grid.}
  \label{fig:Locking}
\end{figure}

Actually the locking effect is related to the pollution of the numerical method precision. This may be explained by a careful analysis of the truncation errors. The solution defined by Eq.~\eqref{eq:J79a} can be recast into
\begin{subequations}\label{sys:def:aniso:corrected}
\begin{equation}
\phi^\eps(x,y) = \phi^0(x,y) + {\varepsilon} \phi_1(x,y) \,, \qquad \nabla_\parallel \phi^0 = 0 \,.
\end{equation}
The property of $\phi^0$ to have no parallel gradients is not exactly met for discretized operators. For instance, the classical finite element spaces do not contain functions with no gradients along $\bb$. Denoting $\Phi^{0,h}$ the discrete approximation of $\phi^0$, the precision of the method is polluted by the discretization error associated to $\eps^{-1}\Delta_\parallel^h \Phi^{0,h}$, $\Delta_\parallel^h$ denoting a discrete approximation of the parallel Laplacian. The truncation error stemming from this discretization is amplified. Therefore, for small values of $\eps$, this term is dominant and large in magnitude compared to the other terms of the left hand side of the discrete equation. Since this contribution is missing in the right hand side of the original discrete equation, the source term $\mathbf{f}^{h}$ being bounded with respect to $\eps$, the norm of the solution is roughly the reciprocal of this amplified discretization error when $\eps \to 0$. This explains the vanishing of the solution norm with $\eps$-values observed on Fig.~\ref{fig:Locking:a}. As outlined in Sec.~\ref{sec:pollution:perp}, refining the mesh may improve the quality of the numerical approximation, the approximation error being proportional to $h^2/\eps$.

To assess further this issue, similar computations are carried out but with a modified equation in which the source term is corrected according to 
\begin{equation}\label{eq:def:aniso:corrected}
\begin{split}
 & -\Delta^h_\perp \Phi^{h} - \frac{1}{\eps}\Delta_\parallel^h \Phi^h = \mathbf{f}^h - \frac{1}{\eps}\Delta_\parallel^h \Phi^{0,h} \,, \\
 &\quad \Phi^{0,h}_{i,j} = \phi^0(x_i,y_j) \,, \quad \mathbf{f}^{h}_{i,j} = f^\eps(x_i,y_j) \,.
 \end{split}
\end{equation}
\end{subequations}
In this corrected equation, the source term is augmented with the discretization error of the parallel laplacian applied to $\phi^0$, which defines a numerical approximation of zero. 
In the corrected equation~\eqref{eq:def:aniso:corrected}, the amplified discretization error of the left hand side is balanced by the correction supplementing the right hand side which restores the accuracy of the approximation regardless of the $\eps$-values, as observed on Fig.~\ref{fig:Locking:b}. This assesses that the locking effect as well as the pollution of the numerical method precision are related to the same amplification (by $\eps^{-1}$) of the approximation error.

Note also that the precision of the numerical approximations issued by a finite difference discretization of the anisotropic equation is deteriorated as soon as $\eps < 1$. Though a convergence may be observed when the mesh is refined, the precision obtained thanks to a discretization of the anisotropic equation is not optimal even for the weakest anisotropies ($\eps=10^{-2}$) reported on Fig.~\ref{fig:Locking:a}. Indeed the  computations carried out with the same discrete operators but a rescaled parallel dynamic offer a significantly improved accuracy.

Note, that for these computations performed with a varying anisotropy direction, the asymmetric scheme \cite{gunter_modelling_2005} is observed to be more accurate than the symmetric scheme (see Rem.~\ref{remark:symetric:scheme}).

\subsubsection{High order discretizations versus rescaled parallel dynamics}

The gain of a higher order discretization is now investigated on computations similar to that of the precedent section. The finite volume framework of \cite{crouseilles_comparison_2015} is investigated for the discretization of the anisotropic problem with second and fourth order discretizations (see Sec.~\ref{sec:Finite:Volume}). 
Following the conclusions of the first numerical experiments, it seems interesting to increase the precision of the discrete parallel Laplacian to reduce the truncation error at the origin of the pollution. To investigate this idea, a mixed precision scheme is implemented. It consists of a second order accurate discretization of the perpendicular Laplacian ($\Delta_\perp^{\VFT}$) while the parallel Laplacian is discretized thanks to a fourth order discretization ($\Delta_\parallel^\VFF$). The numerical approximations issued from this mixed precision scheme are reported on Fig.~\ref{fig:sp:fv_ba2_bp1} together with a second and fourth order discretization of the anisotropic equation.
\def\xxx{0.33\textwidth}
\begin{figure}[!ht]
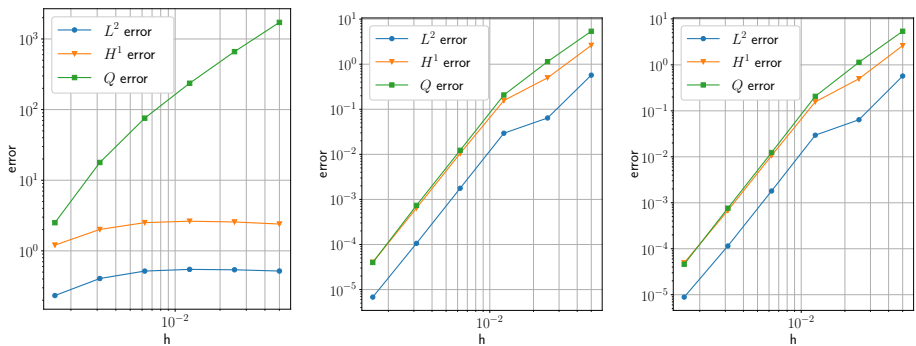
\centering
  \subfloat[Second order, $\varepsilon =10^{-6}$.\label{fig:sp:fv_ba2_bp1:a}] {\resizebox{\xxx}{!}{\input{fig/sp_eps-6_ba2_bp1_deg2.tex}}{}}\hspace*{0.005\textwidth}
  \subfloat[Mixed order (2,4), $\varepsilon =10^{-6}$.\label{fig:sp:fv_ba2_bp1:b}] {\resizebox{\xxx}{!}{\input{fig/sp24_eps-6_ba2_bp1_deg2.tex}}{}}\hspace*{0.005\textwidth}
  \subfloat[Fourth order, $\varepsilon =10^{-6}$.\label{fig:sp:fv_ba2_bp1:c}] {\resizebox{\xxx}{!}{\input{fig/sp_eps-6_ba2_bp1_deg4.tex}}{}}\\

  \caption{$L^2$-norm and $H^1$-norm of the absolute error for the solution approximation $\phi$ together with the $L^2$-norm of the absolute error of the reconstructed flux $\mathcal{Q}^{\eps,h}$ (see Eq.~\eqref{eq:def:reconstructed:flux})
    as functions of the mesh size $h$ for selected values of $\varepsilon$ carried out by second and fourth order finite volume discretizations of the anisotropic problem. The set-up is that of a slowly varying anisotropy direction defined by Eqs.~\eqref{eq:def:setup} with
    $\theta =2, m=1$.}
  \label{fig:sp:fv_ba2_bp1}
\end{figure}

The numerical approximations computed thanks to a second order scheme are not precise enough even on the most refined mesh composed of $640\times 640$ cells (Fig.~\ref{fig:sp:fv_ba2_bp1:a}). The convergence of the numerical method is not observed for the coarsest meshes. This means that the truncation error issued from the discrete parallel Laplacian dominates all the other contributions in the discretized anisotropic equation. This is similar to what is observed on the plots of Fig.~\ref{fig:Locking:a} with no improvements on the error plots when refining the mesh (from $100\times100$ to $300\times300$) for the range of $\eps$-values in $[10^{-9},10^{-6}]$.

The accuracy of the computations carried out thanks to the mixed second and fourth order discretization are noticeably improved. However, the precision of the numerical method is not optimal. Indeed, the solution develops a single oscillation in the perpendicular direction across the computational domain.  The gradients of the $\bb$-field are also very smooth for the selected parameters ($\theta =2, m =1$). The variations of both the function and the $\bb$-field should be well accounted for by a mesh with $30\times30$ cells. Nonetheless, the plots displayed on Fig.~\ref{fig:sp:fv_ba2_bp1:b} show a poor accuracy of the numerical approximations. Here again, this is due to the pollution. Indeed, the convergence rate is observed to be that of a fourth order method, while we could expect that the perpendicular discrete operator drives the precision of the computations with a second order rate of convergence. This proves that the precision of the scheme is limited by the discretization error of the parallel operator.
The plots of Fig.~\ref{fig:sp:fv_ba2_bp1:c} are related to computations performed with a fourth order discretization of both the parallel and perpendicular Laplacians. No gains are observed compared to the computations issued from the mixed precision scheme displayed on Fig.~\ref{fig:sp:fv_ba2_bp1:b}. This confirms that the precision of the approximation is limited by the amplified truncation error originating from the discretization of the parallel Laplacian.

The computations carried out thanks to a rescaled parallel dynamic are displayed on Fig.~\ref{fig:tfi:fv_ba2_bp1}.
The accuracy of the numerical approximation is dramatically improved by two to three order of magnitudes.
For the fourth order discretizations, the solution and its derivatives are approximated with an error smaller than $10^{-3}$ on a $30\times 30$ mesh for computations carried out thanks to a rescaled parallel dynamic (see Fig.~\ref{fig:tfi:fv_ba2_bp1:b}). To gain this accuracy with a discretization of the anisotropic problem (see Figs.~\ref{fig:sp:fv_ba2_bp1:b} and \ref{fig:sp:fv_ba2_bp1:c}) a mesh $300\times 300$ is mandatory. This amounts to a system matrix with a number of rows one hundred times larger and a computational efficiency by no means comparable. This overhead would be even larger for three dimensional computations. The mesh coarsening parameter $g>1$ one can expect thanks to the use of a method with no pollution can be estimated by matching the precision of the non polluted method carried on a coarsened mesh $(h\cdot g)^p$ with that of the polluted scheme $h^p/\eps$, $p$ denoting the approximation order of the method. This yields 
\begin{equation}
g \sim \left({\frac{1}{\eps}}\right)^{1/p}\,.
\end{equation}
The gain is less important when increasing the approximation order but more significant with the increase of the anisotropy strength. For $\eps=10^{-8}$ and $p=4$, the coarsening parameter $g$ is as large as $10^2$. This brings huge savings for the  computational resources, with a number of grid points divided by $10^4$ for two dimensional computations,  to $10^6$ for three dimensional experiments.

We now investigate the approximation of the flux. A comparison of two reconstructions of  $\mathcal{Q}^{\eps}$ can be conducted, with
\begin{subequations}\label{eq:def:fluxes}
\begin{equation}\label{eq:def:reconstructed:flux}
\mathcal{Q}^{\eps,h} = \nabla^h_\perp \Phi^h + \frac{1}{\eps}\nabla^h_\parallel \mathbf{\Phi}^h \,.
\end{equation}
reconstructed thanks to a numerical approximation of the solution $\Phi^h$ and a second approximation using the  auxiliary variable $\mathbf{q}^h$ in order to rescale the parallel dynamics, yielding
\begin{equation}\label{eq:def:rescaled:flux}
\mathcal{Q}^{RPD,h} = \nabla^h_\perp \Phi^h + \nabla^h_\parallel \mathbf{q}^h \,,
\end{equation}
\end{subequations}
The plots of Figs.~\ref{fig:tfi:fv_ba2_bp1} highlight that, whatever the order of the discretizations, the reconstruction with a rescaled parallel dynamic is by several orders of magnitude more accurate than the regular reconstruction $\mathcal{Q}^{\eps,h}$. Actually the plots of $\mathcal{Q}^{\eps,h}$ on Fig.~\ref{fig:tfi:fv_ba2_bp1:b} are very similar to that of Figs.~\ref{fig:sp:fv_ba2_bp1:b} and \ref{fig:sp:fv_ba2_bp1:c}. This puts in perspective that a good approximation of the solution, in terms the $H^1$-norm error does not guaranty an accurate reconstruction of the parallel gradients. This plots emphasize the substantial gains obtained thanks to the rescaling. The flux $\mathcal{Q}^{RPD,h}$ is computed with a precision comparable to the derivatives of the problem solution and is not altered by the anisotropy strength. Contrariwise, the error on $\mathcal{Q}^{\eps,h}$ increases linearly with the values of $1/\eps$ when decreasing $\eps$ from $10^{-6}$ to $10^{-16}$. Note that this latter value of $\eps$ defines anisotropy strengths much more severe than required by the physics, it is therefore a stringent benchmark for the numerical methods. 
\def\xx{0.45\textwidth}
\begin{figure}[!htbp]\centering
  \subfloat[Second order scheme ($\omega=1$).\label{fig:tfi:fv_ba2_bp1:a}] {\resizebox{\xx}{!}{\includegraphics[width=0.5\textwidth]{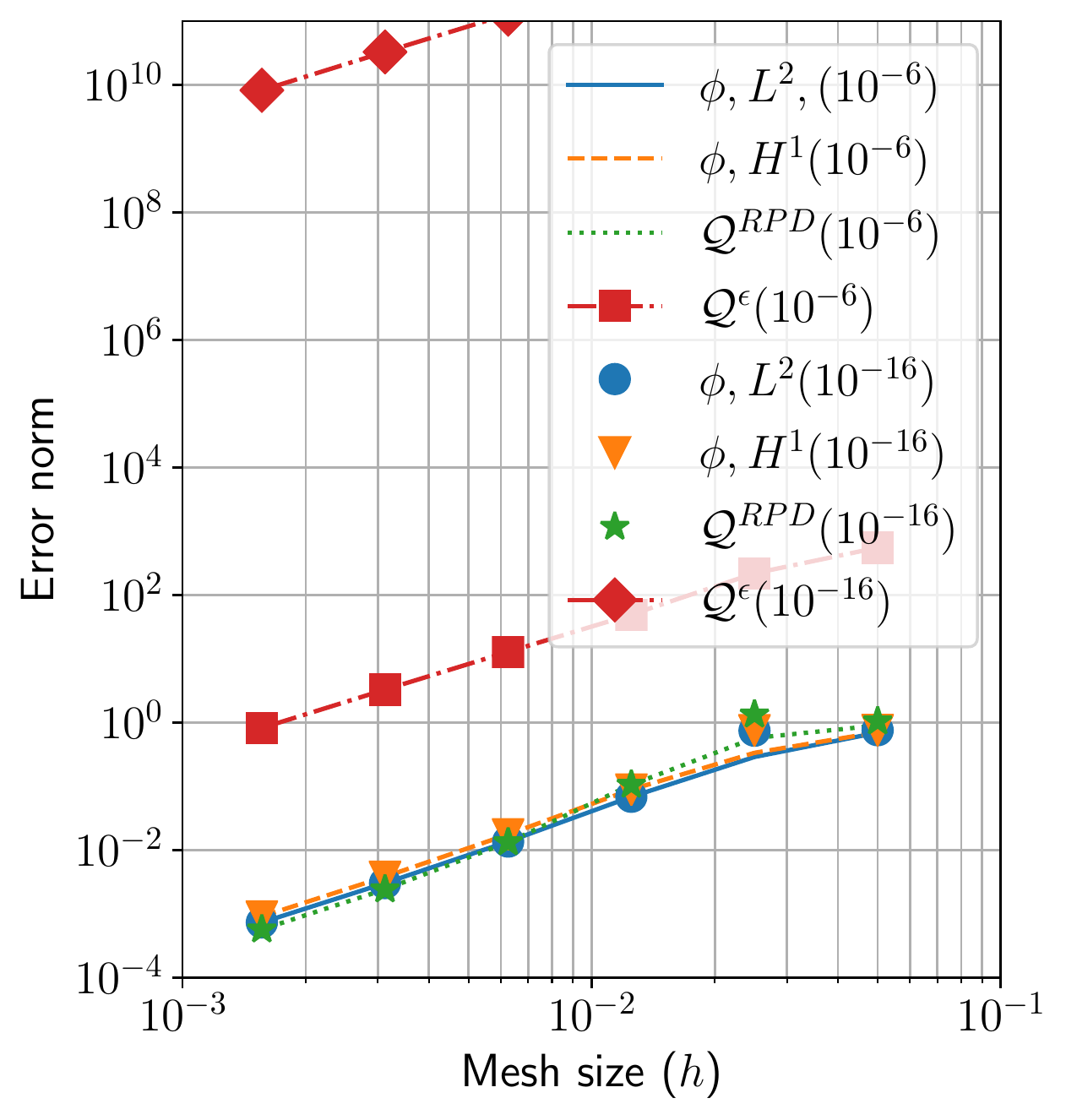}}{}}\hspace*{0.05\textwidth}
  \subfloat[Fourth order scheme ($\omega=1$).\label{fig:tfi:fv_ba2_bp1:b}] {\resizebox{\xx}{!}{\includegraphics[width=0.5\textwidth]{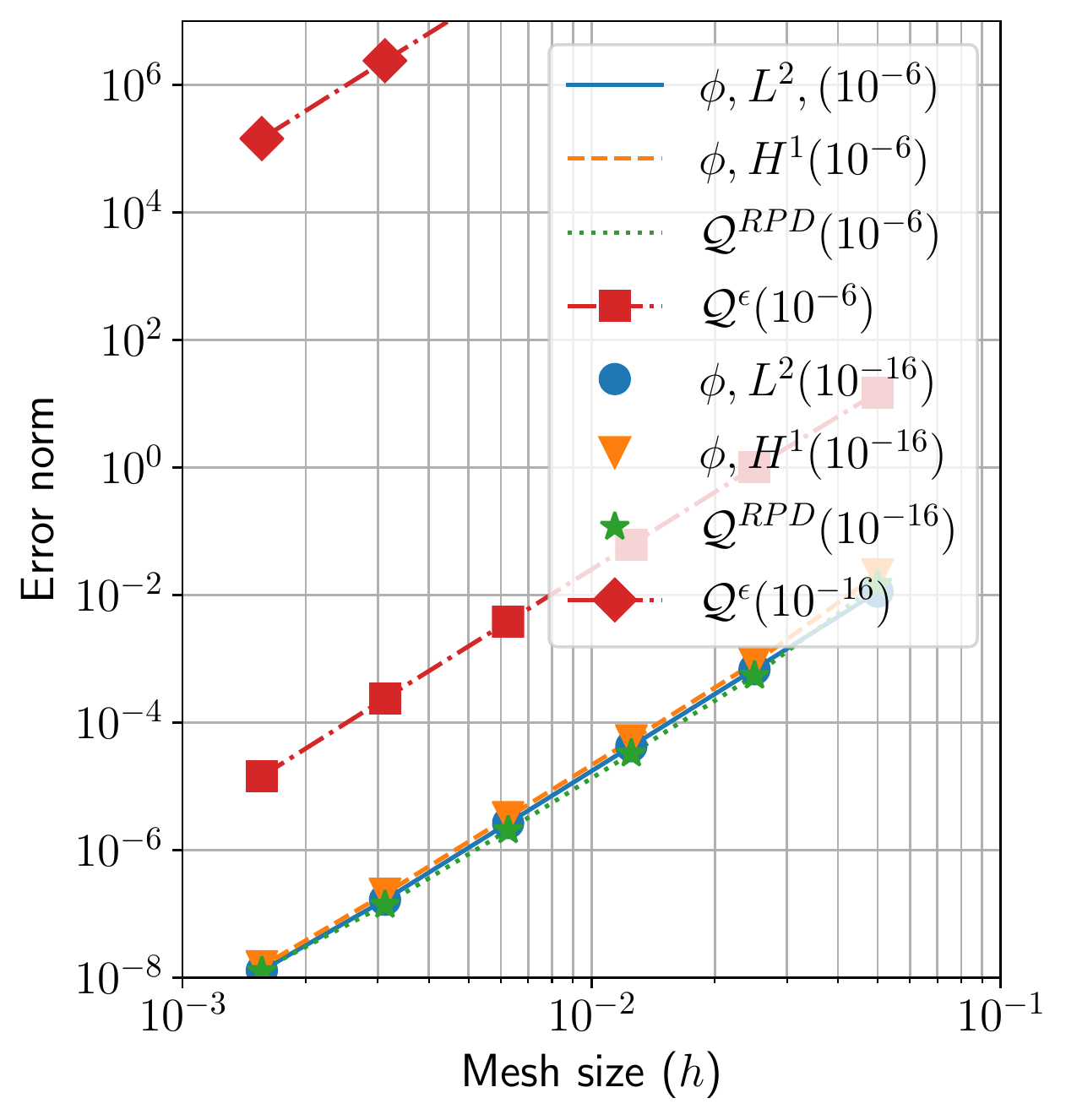}}{}}\\

\subfloat[Anisotropic problem, Fourth order ($\omega=10$).\label{fig:tfi:fv_ba2_bp1:c}] {\resizebox{\xx}{!}{\includegraphics[width=0.5\textwidth]{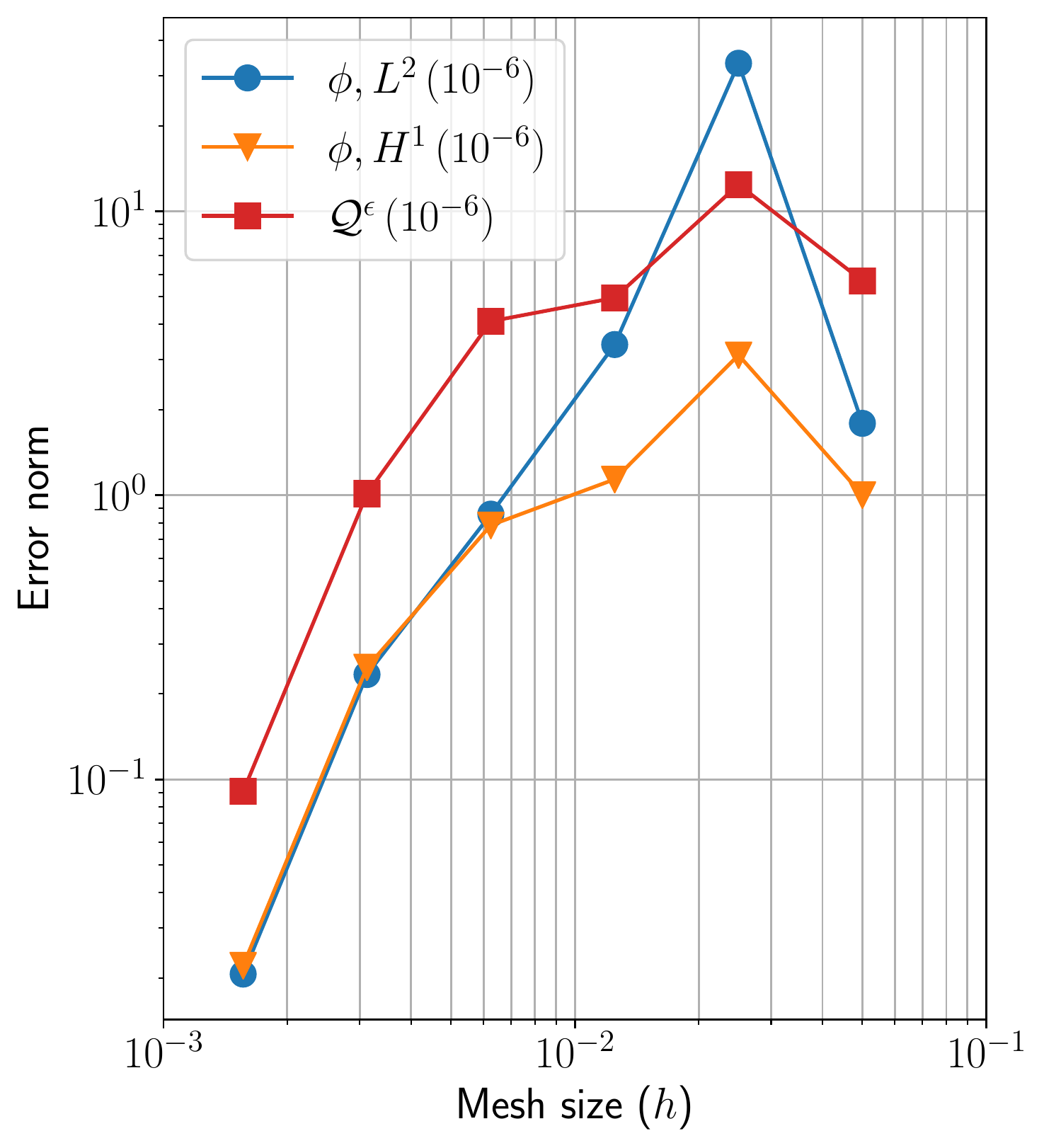}}{}}\hspace*{0.05\textwidth}
  \subfloat[Rescaled Parallel Dynamic, Fourth order ($\omega=10$).\label{fig:tfi:fv_ba2_bp1:d}] {\resizebox{\xx}{!}{\includegraphics[width=0.5\textwidth]{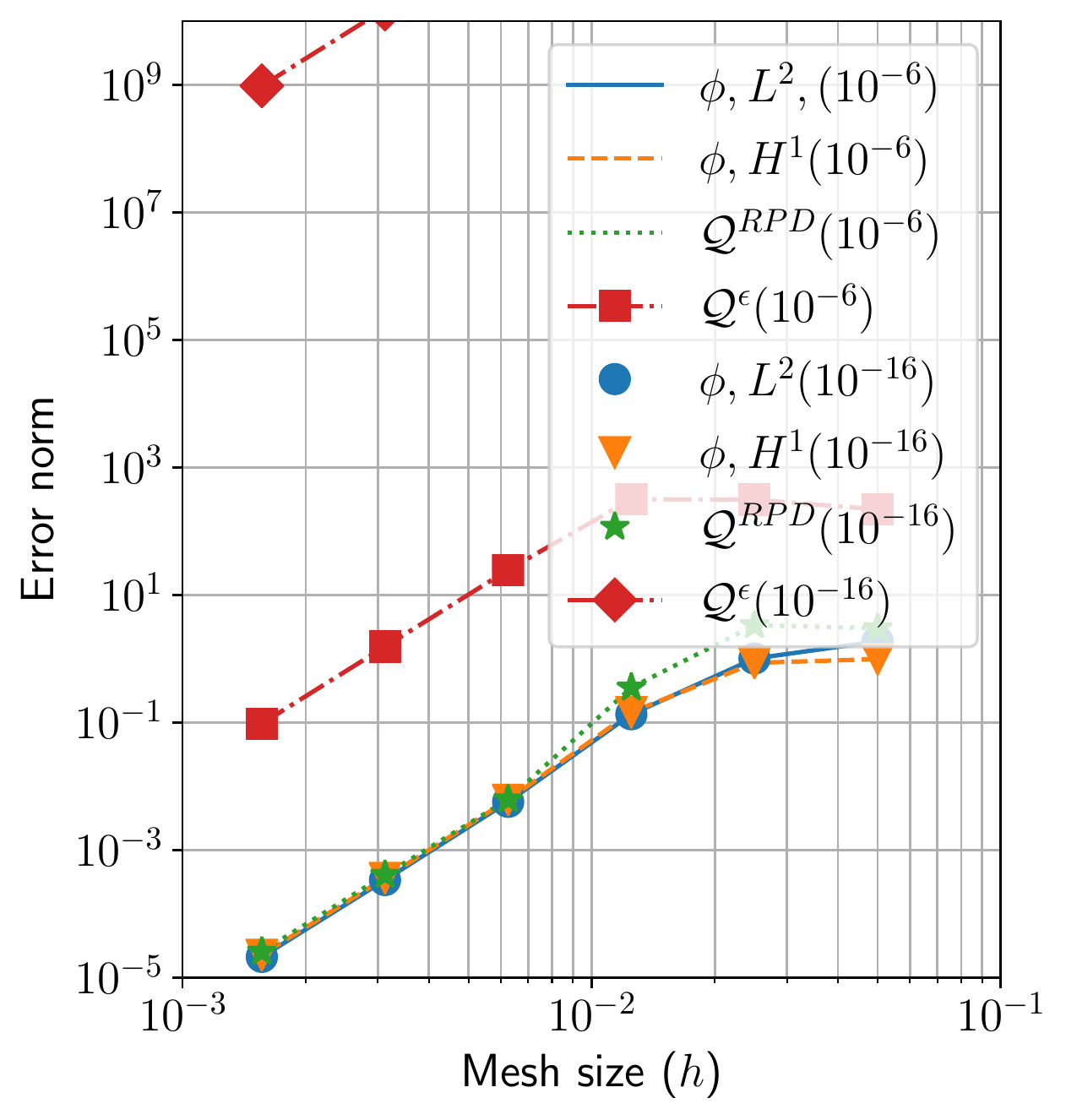}}{}}\\

  \caption{$L^2$-norm and $H^1$-norm of the relative error for the solution approximation $\phi$ together with the $L^2$-norm of the relative error of the reconstructed flux $\mathcal{Q}^{\eps,h}$ and $\mathcal{Q}^{RPD,h}$ (see Eqs.~\eqref{eq:def:fluxes}) as functions of the mesh size $h$ for values of $\varepsilon=10^{-6}$ and $10^{-16}$ carried out by second and fourth order finite volume discretizations of a rescaled parallel dynamic (the Two-Field iterated) system as well as a fourth order discretization of the anisotropic problem. The set-up is that of a slowly varying anisotropy direction defined by Eqs.~\eqref{eq:def:setup} with
    $\theta =2, m=1$ and either $\omega=1$ or $\omega=10$.}
  \label{fig:tfi:fv_ba2_bp1}
\end{figure}

It is important to emphasize that the use of high order methods may not be sufficient. The error plots reported on Figs.~\ref{fig:tfi:fv_ba2_bp1:c} and \ref{fig:tfi:fv_ba2_bp1:d} show a poor accuracy for the approximation of a solution with large perpendicular gradients ($\omega=10$) and moderate anisotropy ($\eps=10^{-6}$). The computations carried without a rescaling of the parallel dynamic (Fig.~\ref{fig:tfi:fv_ba2_bp1:c}) are meaningless. This is totally in line with the conclusions of \cite{crouseilles_comparison_2015} reporting a loss of accuracy for computations with moderate anisotropies ($10^{-6}< \eps < 10^{-3}$). 

The effectiveness of high order methods is investigated further with Finite Element discretizations implemented on the same set-up with $\mathbb{Q}_2$, $\mathbb{Q}_4$ and $\mathbb{Q}_6$ methods, defining third, fifth and seventh order discretizations of the solution in the $L^2$-norm. For the $H^1$-norm of the solution error and the $L^2$-norm of the fluxes error these approximation are one order of magnitude less precise.  The plots are gathered in Fig.~\ref{fig:fe_ba2_bp1} for $\eps$-values equal to $10^{-6}$ and $10^{-16}$. 
\def\xxx{0.32\textwidth}
\begin{figure}[!htbp]

  \subfloat[$\mathbb Q_2$, $\varepsilon =10^{-6}$] {\resizebox{\xxx}{!}{\input{fig/mm_eps1e-06_ba2_bp1_deg2.tex}}{}}\hspace*{0.01\textwidth}\hspace*{0.005\textwidth}
  \subfloat[$\mathbb Q_4$, $\varepsilon =10^{-6}$\label{fig:fe_ba2_bp1:b}] {\resizebox{\xxx}{!}{\input{fig/mm_eps1e-06_ba2_bp1_deg4.tex}}{}}\hspace*{0.005\textwidth}
  \subfloat[$\mathbb Q_6$, $\varepsilon =10^{-6}$\label{fig:fe_ba2_bp1:c}] {\resizebox{\xxx}{!}{\input{fig/mm_eps1e-06_ba2_bp1_deg6.tex}}{}}\\

  \subfloat[$\mathbb Q_2$, $\varepsilon =10^{-16}$]{\resizebox{\xxx}{!}{\input{fig/mm_eps1e-16_ba2_bp1_deg2.tex}}{}}\hspace*{0.01\textwidth}\hspace*{0.005\textwidth}
  \subfloat[$\mathbb Q_4$, $\varepsilon =10^{-16}$]{\resizebox{\xxx}{!}{\input{fig/mm_eps1e-16_ba2_bp1_deg4.tex}}{}}\hspace*{0.01\textwidth}\hspace*{0.005\textwidth}
  \subfloat[$\mathbb Q_6$, $\varepsilon =10^{-16}$]{\resizebox{\xxx}{!}{\input{fig/mm_eps1e-16_ba2_bp1_deg6.tex}}{}}\\

  \caption{$L^2$-norm and $H^1$-norm of the absolute error for the solution approximation $\phi$ together with the $L^2$-norm of the absolute error of the reconstructed flux $\mathcal{Q}^{\eps,h}$ and $\mathcal{Q}^{RPD,h}$ (see Eqs.~\ref{eq:def:fluxes}) as functions of the mesh size $h$ for values of $\varepsilon=10^{-6}$ and $10^{-16}$ carried out by $\mathbb Q_2$, $\mathbb Q_4$ and $\mathbb Q_6$ Finite Element discretizations of a rescaled parallel dynamic (Micro-Macro) system. The set-up is that of a slowly varying anisotropy direction defined by Eqs.~\eqref{eq:def:setup} with
    $\theta =2, m=1$.}
  \label{fig:fe_ba2_bp1}
\end{figure}
Increasing the precision of the numerical method thanks to high order discretizations improves the quality of the approximations. Indeed $\mathbb{Q}_4$ and $\mathbb{Q}_6$ finite element computations, as shown on Figs.~\ref{fig:fe_ba2_bp1:b} and \ref{fig:fe_ba2_bp1:c}, permit the reconstruction of quite an accurate flux approximation $\mathcal{Q}^{\eps,h}$. This is true on the most refined meshes for both orders, however with a noticeably deteriorated precision compared to the rescaled flux approximation $\mathcal{Q}^{RPD,h}$. Moreover, these computations do not show the global picture, since the precision of 
$\mathcal{Q}^{\eps,h}$ deteriorates with either vanishing $\eps$ and the regularity of the solution, or more specifically with the magnitude of its perpendicular derivatives as shown on Fig.~\ref{fig:tfi:fv_ba2_bp1}. The use of an oversampled mesh is mandatory to carry out an approximation with a sufficient accuracy at the price of the computational cost due to the increased size of the system matrix together with its condition number.  

\begin{figure}[!ht]
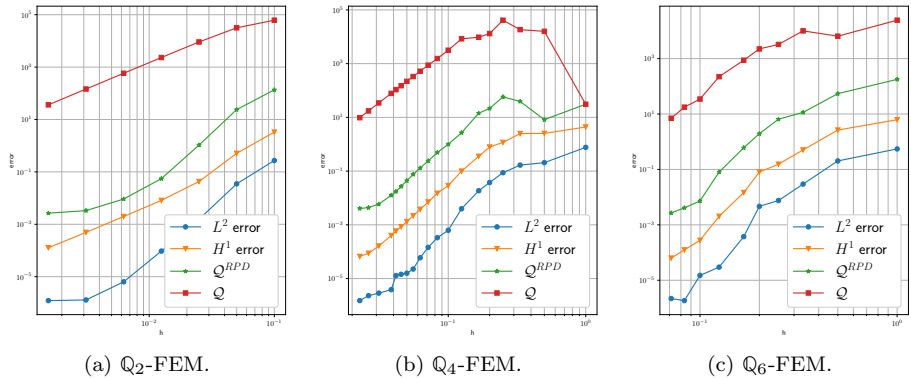
\centering
\def\xxx{0.32\textwidth}
  \subfloat[$\mathbb Q_2$-FEM.] {\resizebox{\xxx}{!}{\input{fig/tfi_eps0_0.001_it10_eps1e-06_ba10_bp2_deg2.tex}}{}}\hspace*{0.01\textwidth}\hspace*{0.005\textwidth}
  \subfloat[$\mathbb Q_4$-FEM.] {\resizebox{\xxx}{!}{\input{fig/tfi_eps0_0.001_it10_eps1e-06_ba10_bp2_deg4.tex}}{}}\hspace*{0.01\textwidth}\hspace*{0.005\textwidth}
  \subfloat[$\mathbb Q_6$-FEM.] {\resizebox{\xxx}{!}{\input{fig/tfi_eps0_0.001_it10_eps1e-06_ba10_bp2_deg6.tex}}{}}\\

  \caption{$L^2$-norm and $H^1$-norm of the absolute error for the solution approximation $\phi$ together with the $L^2$-norm of the absolute error of the reconstructed flux $\mathcal{Q}^{\eps,h}$ and $\mathcal{Q}^{RPD,h}$ (see Eqs.~\ref{eq:def:fluxes}) as functions of the mesh size $h$ for  $\varepsilon=10^{-6}$ carried out by $\mathbb Q_2$, $\mathbb Q_4$ and $\mathbb Q_6$ Finite Element discretizations of a rescaled parallel dynamic (TFI) system. The set-up is defined by Eqs.~\eqref{eq:def:setup} with
    $\theta =10, m=2$ including close field lines.}
  \label{fig:fe_ba10_bp2}
\end{figure}

With the last set-up the efficiency of numerical methods is experienced in the presence of closed field lines. This is selected by setting $\theta =10$, $m=2$ into Eqs.~\eqref{eq:def:setup}. The field
$B$ is equal to $(0,0)^T$ in the middle of the domain located at
$x=1/2$, $y=1/2+\pi/20$. One difficulty of this problem relies therefore in the presence of a point where the field $\bb$ is not defined. The workaround implemented for the computations consists in setting $b=(0,1)^T$ in this point.
The computations related to this set-up are reported in Fig.~\ref{fig:fe_ba10_bp2}.
It is challenging to recover the optimal convergence rate for this problem because the gradients of the $\bb$-field captured by the mesh stiffen with the grid resolution. The reconstruction of the flux from the derivatives of the solution approximation ($\mathcal{Q}^{\eps,h}$) cannot produce a meaningful  estimate whatever the approximation order and the mesh size used. For these computations too, the precision of this quantity is affected by an error proportional to reciprocal of
$\eps$. Contrariwise, the flux $\mathcal{Q}^{RPD}$ reconstructed thanks to a rescaled parallel dynamics provides an accurate approximation when the mesh is refined enough.

\def\xx{0.4\textwidth}

Let us now focus on the propagation of round-off errors. To this end,
the framework investigated is that of a coordinate aligned anisotropy which amounts to choose
$\theta=0$ in Eqs.~\eqref{eq:def:setup}.
The most fundamental difference with this setting is that the finite element space contains non trivial functions that are
constant in the direction of the anisotropy: the issue related to the precision pollution 
is therefore expelled from the problem. The computations displayed on Fig.~\ref{fig:fe_ba0_bp1} are carried out thanks to the TFI method implementing a rescaled parallel dynamic. Two series of plots are proposed, the first ones related to $\eps=10^{-6}$ the second ones to a severe anisotropy with $\eps=10^{-16}$. The accuracy of the solution approximation as well as the flux reconstructed with a rescaled parallel dynamic $\mathcal{Q}^{RPD,h}$ are unaffected by the anisotropy strength. The error plots follow the expected convergence rate until a critical mesh size $h_c$ for which the amplified round-off errors match the precision of the discretizations (see Sec.~\ref{sec:arithmetic}). It is important to note that, the value of $h_c$ does not depend on $\eps$: the same threshold is observed on Figs.~\ref{fig:fe_ba0_bp1:a} and \ref{fig:fe_ba0_bp1:e} implementing a $\mathbb{Q}_2$-FEM discretization but different $\eps$-values. The same conclusion holds true for the plots of Figs.~\ref{fig:fe_ba0_bp1:b} and \ref{fig:fe_ba0_bp1:f}  as well as Figs.~\ref{fig:fe_ba0_bp1:c} and \ref{fig:fe_ba0_bp1:f} with the exact same error plots whatever the anisotropy strength. This is an outcome of matrices issued from the discretization of systems with a rescaled parallel dynamic providing an amplification of round-off errors dependant of the mesh size $h$ but unrelated to $\eps$.  Contrariwise, the quantity computed with non rescaled parallel dynamics exhibit a precision dependent on $\eps$ with an amplification of round-off errors proportional to $1/\eps$. This is manifest on the plots related to the error of the reconstructed flux $\mathcal{Q}^{\eps,h}$.

\begin{figure}[!ht]
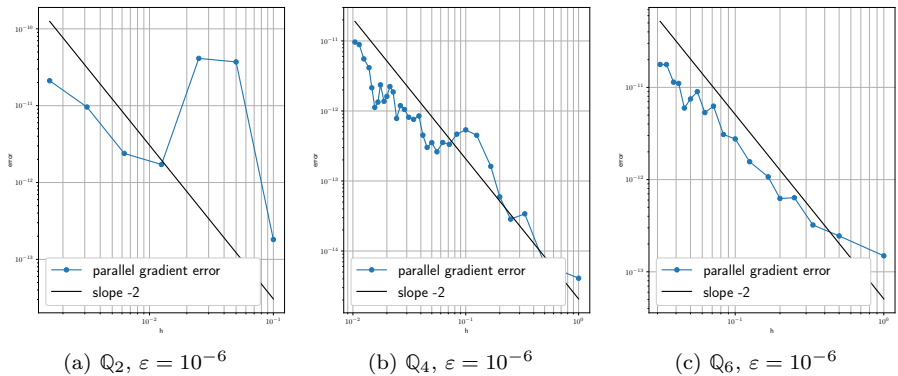
\centering
\def\xxx{0.32\textwidth}
  \subfloat[$\mathbb Q_2$, $\varepsilon =10^{-6}$] {\resizebox{\xxx}{!}{\input{fig/tfi_eps0_0.001_it10_eps1e-16_ba0_bp1_deg2_para.tex}}{}}\hspace*{0.005\textwidth}
  \subfloat[$\mathbb Q_4$, $\varepsilon =10^{-6}$] {\resizebox{\xxx}{!}{\input{fig/tfi_eps0_0.001_it10_eps1e-16_ba0_bp1_deg4_para.tex}}{}}\hspace*{0.005\textwidth}
  \subfloat[$\mathbb Q_6$, $\varepsilon =10^{-6}$] {\resizebox{\xxx}{!}{\input{fig/tfi_eps0_0.001_it10_eps1e-16_ba0_bp1_deg6_para.tex}}{}}\\

  \caption{Absolute errors of the parallel gradient $||\nabla_\parallel \phi
    ^\varepsilon - \nabla_\parallel^h \phi ^\varepsilon _h||_{L^2}$ as functions of mesh size for
    $\varepsilon = 10^{-16}$ for an anisotropy direction aligned with the coordinate system ($\theta =0, m=1$).}
  \label{fig:fe_ba0_bp1_para}
\end{figure}

This feature may be analysed thanks to the investigation of the parallel gradient reconstruction.
The parallel gradient error as a function of the mesh size is plotted on Fig.~\ref{fig:fe_ba0_bp1_para} for $\varepsilon =10^{-16}$. First, we note that the value of $\eps$ is too small for the parallel gradient to be reconstructed at the right scale ($\eps$) whatever the precision order of the discretization. Second, the best precision is achieved for coarsest meshes and scales as the (TFI) matrix
condition number, {\it i.e.} as $1/h ^2$. For instance, 
the $\mathbb{Q}_6$-scheme  achieves its best precision on a
mesh with a single element ($6\times 6 $ discretization points) yielding an absolute error of roughly
$10^{-13}$. This error is amplified by $1/\eps$ when introduced in the reconstructed flux and explains the loss of precision for the approximation of this quantity on refined meshes.
 This suggests that the flux
calculated directly from $\phi ^{\varepsilon,h} $ could be meaningful for
medium-sized meshes with $\varepsilon $ of the order of $10^{-6}$ and outlines the advantage to work a system with a rescaled parallel dynamic.

\begin{figure}[!ht]\centering
\def\xxx{0.32\textwidth}
  \subfloat[$\mathbb Q_2$, $\varepsilon =10^{-6}$\label{fig:fe_ba0_bp1:a}] {\resizebox{\xxx}{!}{\input{fig/tfi_eps0_0.001_it10_eps1e-06_ba0_bp1_deg2.tex}}{}}\hspace*{0.005\textwidth}
  \subfloat[$\mathbb Q_4$, $\varepsilon =10^{-6}$\label{fig:fe_ba0_bp1:b}] {\resizebox{\xxx}{!}{\input{fig/tfi_eps0_0.001_it10_eps1e-06_ba0_bp1_deg4.tex}}{}}\hspace*{0.005\textwidth}
  \subfloat[$\mathbb Q_6$, $\varepsilon =10^{-6}$\label{fig:fe_ba0_bp1:c}] {\resizebox{\xxx}{!}{\input{fig/tfi_eps0_0.001_it10_eps1e-06_ba0_bp1_deg6.tex}}{}}\\

  \subfloat[$\mathbb Q_2$, $\varepsilon =10^{-16}$\label{fig:fe_ba0_bp1:e}]{\resizebox{\xxx}{!}{\input{fig/tfi_eps0_0.001_it10_eps1e-16_ba0_bp1_deg2.tex}}{}}\hspace*{0.005\textwidth}
  \subfloat[$\mathbb Q_4$, $\varepsilon =10^{-16}$\label{fig:fe_ba0_bp1:f}]{\resizebox{\xxx}{!}{\input{fig/tfi_eps0_0.001_it10_eps1e-16_ba0_bp1_deg4.tex}}{}}\hspace*{0.005\textwidth}
  \subfloat[$\mathbb Q_6$, $\varepsilon =10^{-16}$\label{fig:fe_ba0_bp1:g}]{\resizebox{\xxx}{!}{\input{fig/tfi_eps0_0.001_it10_eps1e-16_ba0_bp1_deg6.tex}}{}}\\

  \caption{$L^2$-norm and $H^1$-norm of the absolute error of the solution approximation $\phi$ together with the $L^2$-norm of the absolute error of the reconstructed flux $\mathcal{Q}^{\eps,h}$ and $\mathcal{Q}^{RPD,h}$ (see Eqs.~\ref{eq:def:fluxes}) as functions of the mesh size $h$ for  either $\eps=10^{-6}$ or $\eps=10^{-16}$ carried out by $\mathbb Q_2$, $\mathbb Q_4$ and $\mathbb Q_6$ Finite Element discretizations of a rescaled parallel dynamic (TFI) system. The set-up is defined by Eqs.~\eqref{eq:def:setup} with
    $\theta =0, m=1$. \label{fig:fe_ba0_bp1}}
\end{figure}

\section{Conclusions}
In this paper, the precision pollution arising in the numerical approximation of anisotropic problems is investigated. 
The origin of this loss of accuracy is related to an amplification of the truncation error of the parallel operator. The analysis carried out in this paper demonstrates that the use of high order schemes, classically operated to offset this error amplification, does not overcome this difficulty. Indeed, though the convergence of these methods may be observed at the expected rate, their precision remains far from optimal because it is altered by an error amplified by the imbalance between the perpendicular and parallel diffusions. The numerical investigations conducted within this document show a discrepancy of the accuracy of numerical approximations, by many order of magnitudes, depending on the strength of the anisotropy. This prevents from computing meaningful numerical approximations for severe anisotropies unless the mesh is oversampled in order to alleviate the amplification of the approximation error. 

The most stringent criterion for validating the effectiveness of a numerical method reveals to be the error $L^2$-norm of the flux associated to the anisotropic problem. In this quantity, the parallel gradients are rescaled by the coefficient defining the anisotropy strength ($\eps^{-1}\gg 1$). Therefore, the approximation of this quantity requires a precise reconstruction of the solution parallel gradients which is the most challenging difficulty raised by these anisotropic problems. 

We introduce a new approach consisting in rescaling the parallel dynamic in order to remove the stiffness from the equation and cancel the cause of the error amplification. This is achieved thanks to an auxiliary variable accounting for the parallel gradients. The numerical experiments performed on numerous and demanding benchmarks demonstrate the effectiveness of the approaches proposed within this paper, with numerical approximations of the anisotropic problem solution, as well as the reconstructed flux, unaltered by the anisotropy strength. This makes possible the computations on meshes refined according to the gradients of the solution rather than the need to alleviate the amplification of the approximation error.

%

\section*{Acknowledgements}
This work has been supported by the french ``Agence Nationale pour la
  Recherche (ANR)'' in the frame of the contract ANR-11-MONU-009-01 ``MOONRISE: MOdels, Oscillations and NumeRIcal SchEmes'' (2015-2019) as well as the ``labex CIMI'' (International Centre for Mathematics and Computer Science in Toulouse) in the frame of the project ``SCANISO: SCalable solvers for ANISOtropic equations arising in magnetized plasma simulations'' (2017-2019).
Support from the ``F\'ed\'eration de Fusion pour la Recherche par Confinement Magn\'etique'' (FrFCM) in the frame of the project ``NEMESIA: Numerical mEthods for 
Macroscopic models of magnEtized plaSmas and related anIsotropic equAtions''.\newline
CY is supported by the Fog Research Institute under contract no. FRI-454.
CY is also supported by Heilongjiang Natural Science Foundation (LH2019A013).\newline
FD acknowledges invitations as a scholar professor by Harbin Institute of Technology in 2016, 2017, 2018 and 2019.
CY acknowledges invitations as a scholar professor by Universit\'e Paul Sabatier, Toulouse 3, during 2019.
\newline The authors acknowledge fruitful discussions with Alexei Lozinski.
\bibliographystyle{abbrv}
\bibliography{bib}

%
  \appendix

\section{Discrete differential operators}
\subsection{Mesh definition, notations}\label{sec:appendix:notations}
The computational domain $\Omega=(0,1)\times (0,1)$ is decomposed into cells 
\begin{equation*}
K_{i,j}  =  [x_{i-1/2},x_{i+1/2}] \times [y_{j-1/2},y_{j+1/2}], \quad (i,j)\in\{1,\dots,N\} \times \{1,\dots,N\},
\end{equation*}
where 
\begin{eqnarray*}
x_{i+1/2} & = & (i+1/2)h,\quad i\in\{0,\dots,N\}, \\[3mm]
y_{j+1/2} & = & (j+1/2)h,\quad j\in\{0,\dots,N\}, 
\end{eqnarray*}
 the mesh sizes being defined as
\begin{equation}
h =  \frac{1}{N}\,.
\end{equation}
The faces of the control volumes are 
\begin{eqnarray*}
\sigma^x_{i} & = & [x_{i-1/2},x_{i+1/2}],\quad i\in\{1,\dots,N\}, \\
\sigma^y_{j} & = & [y_{j-1/2},y_{j+1/2}],\quad k\in\{1,\dots,N\}.
\end{eqnarray*}

We note $\Phi^h$ the vector of the function point values at the cell centers and $\bar \Phi^h$ the vector of the cell averages with
\begin{equation}
  \left( \Phi^h\right)_{i,j} = \phi_{i,j} \approx \phi(x_i,y_j) \,, \qquad
  \left( \bar \Phi^h\right)_{i,j} = \bar\phi_{i,j} \approx \frac{1}{h^2} \int_{K_{i,j}} \phi(x,y) dx dy \,
\end{equation}
and
\begin{align}
  A_\parallel =  b \otimes b  \,, \qquad A_\perp =  \Id - A_\parallel
\end{align}
\subsection{Finite Difference discretizations}\label{sec:finite:difference}
Two finite difference discretizations are considered within this work. They are derived from the so-called symmetric and anti-symmetric fluxes introduced in \cite{gunter_modelling_2005}. 
The anti-symmetric schemes relies on a definition of the flux components at the cell interfaces with
\begin{subequations}\label{eqs:Asym}
\begin{equation}
\begin{split}
  \left(\Delta_\parallel^\ASYM \Phi^h\right)_{i,j}&= \frac{1}{h} \Big( (F_{\parallel,x}^\ASYM)_{i+1/2,j} -  (F_{\parallel,x}^\ASYM)_{i-1/2,j} \\&\qquad + (F_{\parallel,y}^\ASYM)_{i,j+1/2} -  (F_{\parallel,y}^\ASYM)_{i,j-1/2}\Big) \,,
\end{split}
\end{equation}

where 
\begin{align}
\begin{split}
(F_{\parallel,x}^\ASYM)_{i+1/2,j} =&\left(A_{\parallel,xx}\right)_{i+1/2,j} \left(\partial_x^\ASYM \Phi^h\right)_{i+1/2,j} +\\&\qquad \left(A_{\parallel,xy}\right)_{i+1/2,j} \left(\partial_y^\ASYM \Phi^h\right)_{i+1/2,j} \,,
\end{split}\\
\begin{split}
(F_{\parallel,y}^\ASYM)_{i,j+1/2} =&\left(A_{\parallel,yx}\right)_{i,j+1/2} \left(\partial_x^\ASYM \Phi^h\right)_{i,j+1/2} +\\& \qquad \left(A_{\parallel,yy}\right)_{i,j+1/2} \left(\partial_y^\ASYM \Phi^h\right)_{i,j+1/2} \,,
\end{split}
\end{align}
and
\begin{equation}
\begin{split}
  \left(\partial_x^\ASYM \Phi^h\right)_{i+1/2,j} &= \frac{1}{h}\left(\phi_{i+1,j} - \phi_{i,j}  \right)\,,\\  \left(\partial_y^\ASYM \Phi^h\right)_{i,j+1/2} &= \frac{1}{h}\left(\phi_{i,j+1} - \phi_{i,j}  \right) \,.
  \end{split}
\end{equation}
\begin{align}
  \left(\partial_x^\ASYM \Phi^h\right)_{i,j+1/2} &= \frac{1}{4h}\Big(\left(\phi_{i+1,j}+\phi_{i+1,j}\right) - \left( \phi_{i-1,j+1}+\phi_{i-1,j}  \right) \Big) \,,\\
  \left(\partial_y^\ASYM \Phi^h\right)_{i+1/2,j} &= \frac{1}{4h}\Big(\left(\phi_{i+1,j+1} + \phi_{i,j+1}\right) - \left( \phi_{i+1,j-1} + \phi_{i,j-1} \right) \Big)\,.
\end{align}
\end{subequations}
For the symmetric scheme, all the flux components are carried out at the same location, yielding:
\begin{subequations}\label{eqs:Sym}
\begin{equation}
\begin{split}
  \left(\Delta_\parallel^\SYM \Phi^h\right)_{i,j}&= \frac{1}{2 h} \Big( (F_{\parallel,x}^\SYM)_{i+1/2,j+1/2} + (F_{\parallel,x}^\SYM)_{i+1/2,j-1/2} \\
   & \qquad - \left( (F_{\parallel,x}^\SYM)_{i-1/2,j+1/2} + (F_{\parallel,x}^\SYM)_{i-1/2,j-1/2} \right)  \\
   &  + (F_{\parallel,y}^\SYM)_{i+1/2,j+1/2}+(F_{\parallel,y}^\SYM)_{i-1/2,j+1/2} \\
   & \qquad - \left((F_{\parallel,y}^\SYM)_{i+1/2,j-1/2}+(F_{\parallel,y}^\SYM)_{i-1/2,j-1/2}  \right)\Big) \,,
\end{split}
\end{equation}

with
\begin{align}
  \left(\partial_x^\SYM \Phi^h\right)_{i+1/2,j+1/2} &= \frac{1}{2h}\Big(\left(\phi_{i+1,j+1}+\phi_{i+1,j}\right) - \left( \phi_{i,j+1}+\phi_{i,j}  \right) \Big) \,,\\
  \left(\partial_y^\SYM \Phi^h\right)_{i+1/2,j+1/2} &= \frac{1}{2h}\Big(\left(\phi_{i+1,j+1} + \phi_{i,j+1}\right) - \left( \phi_{i+1,j} + \phi_{i,j} \right) \Big)\,.
\end{align}
\end{subequations}

The discretization of the perpendicular operator is deduced from that of the parallel Laplacian.

\subsection{Finite Volume discretizations}\label{sec:Finite:Volume}

The Finite Volume discretization is similar to the ones implemented in \cite{crouseilles_comparison_2015} derived from \cite{zhang_fourth-order_2012}. The discrete Laplace operators are obtained thanks to an integration over the control volume $K_{i,j}$ with
\begin{eqnarray*}
\Delta^h_{\bot}\phi_{i,j}&=&\frac{1}{h} \Big({F^h_{\bot,x}(\phi)_{i+1/2,j} - F^h_{\bot,x}(\phi)_{i-1/2,j}}{\Delta x}\Big) \\
&& \hspace*{10em} + \frac{1}{h}\Big({F^h_{\bot,y}(\phi)_{i,j+1/2} - F^h_{\bot,y}(\phi)_{i,j-1/2}} \Big),\\[3mm]
\Delta^h_{\|}\phi_{i,j}&=&\frac{1}{h} \Big(F^h_{\|,x}(\phi)_{i+1/2,j} - F^h_{\|,x}(\phi)_{i-1/2,j}\Big) \\
&& \hspace*{10em} + \frac{1}{h}\Big({F^h_{\|,y}(\phi)_{i,j+1/2} - F^h_{\|,x}(\phi)_{i,j-1/2}}\Big).
\end{eqnarray*}
where the fluxes are defined as the integration along the edges of the control volume  $K_{i,j}$:

\begin{eqnarray*}
F_{\|,x}(\phi)_{i+1/2,j} & = & \frac{1}{h}\int_{\sigma_j^y} (A_{\|,xx} \partial_x\phi) (x_{i+1/2}, y) + (A_{\|,xy} \partial_y\phi) (x_{i+1/2}, y) dy ,\\[3mm]
F_{\|,y}(\phi)_{i,j+1/2} & = & \frac{1}{h}\int_{\sigma_i^x} (A_{\|,xy} \partial_x\phi) (x, y_{j+1/2}) + (A_{\|,yy} \partial_y\phi) (x, y_{j+1/2}) dx.
\end{eqnarray*}
An approximation of $\frac{1}{h} \int_{\sigma_j^y} (A\,\partial_x\phi)(x_{i+1/2},y)dy$ denoted $(A\,\partial_x {\phi})_{i,j+1/2}$ is then introduced. A fourth order as well as a second order approximation of this integrals are used to define the discrete operators. For succintness, only the fourth order  approximation is precised (see \cite{crouseilles_comparison_2015,zhang_fourth-order_2012} for details) yielding the following definition
\begin{multline*}
( A\, \partial_x\phi)_{i+1/2,j} =  A_{i+1/2,j} \partial_x\phi_{i+1/2,j} \\
+ \frac{1}{48}\left(  A_{i+1/2,j+1}  -  A_{i+1/2,j-1} \right) 
 \left( \partial_x\phi_{i+1/2,j+1}  - \partial_x\phi_{i+1/2,j-1} \right) .
\end{multline*}
and
\begin{equation*}
\partial_x\phi_{i+1/2,j} = \frac{5}{4h}(\bar{\phi}_{i+1,j} - \bar{\phi}_{i,j}) - \frac{1}{12h}(\bar {\phi}_{i+2,j} - \bar {\phi}_{i-1,j}).
\end{equation*}

To impose the boundary conditions, fifth order extrapolation formulae are implemented to define the values of two levels of ghost cells. For Dirichlet boundary conditions, the following relations are used to compute the values carried by the ghost cells:
\begin{eqnarray*}
\bar \phi_{0,j} & = & \frac{1}{10}(-87\bar \phi_{1,j} + 63\bar \phi_{2,j} - 37 \bar \phi_{3,j} + 13\bar \phi_{4,j} - 2\bar \phi_{5,j} + 60 \bar \phi_{1/2,j}),\\[3mm]
\bar \phi_{-1,j} & = & \frac{1}{5}(-336\bar \phi_{1,j} + 289\bar \phi_{2,j} - 186\bar \phi_{3,j} + 69\bar \phi_{4,j} - 11\bar \phi_{5,j}+ 175\bar \phi_{1/2,j}).
\end{eqnarray*}

\section{Discrete operators and truncation errors.}\label{Appendix:Tronc}

The computations related in this section are performed under the assumption of a uniform magnetic field defined by 
\begin{equation}
b(x,y) = \left( \begin{array}[c]{c}
\alpha=\cos(a)\\ \beta=\sin(a)
\end{array}\right)
\end{equation} 
First, the toncature error associated to the finite differenced operators defined in Sec. are stated with
\begin{subequations}
\begin{equation}
\begin{split}
\left(\Delta_\parallel^\ASYM \bar \Phi^h\right)_{i,j} &= \Delta_\parallel \phi(x_i,y_j) \\
	&+ \alpha\beta\frac{h^2}{3}\left(\frac{\partial^4}{\partial x^3 \partial y} \phi(x_i,y_j) + \frac{\partial^4}{\partial x \partial y^3} \phi(x_i,y_j)\right) \\
  &+ \frac{h^2}{12}\left(\alpha^2 \frac{\partial^4}{\partial x^4} \phi(x_i,y_j) + \beta^2 \frac{\partial^4}{\partial y^4} \phi(x_i,y_j)\right)+ \mathcal{O}(h^4)\,,
  \end{split}
\end{equation}
and
\begin{equation}
\begin{split}
\left(\Delta_\parallel^\SYM \bar \Phi^h\right)_{i,j} &= \Delta_\parallel \phi(x_i,y_j) \\
&+\alpha^2\frac{h^2}{12}\left( \frac{\partial^4}{\partial x^4} \phi(x_i,y_j)+3 \frac{\partial^4}{\partial x^2\partial y^2} \phi(x_i,y_j)\right)\\
	&+ \alpha\beta\frac{h^2}{3}\left(\frac{\partial^4}{\partial x^3 \partial y} \phi(x_i,y_j) + \frac{\partial^4}{\partial x \partial y^3} \phi(x_i,y_j)\right) \\
	&+\beta^2\frac{h^2}{12}\left( \frac{\partial^4}{\partial y^4} \phi(x_i,y_j)+3 \frac{\partial^4}{\partial x^2\partial y^2} \phi(x_i,y_j)\right)+ \mathcal{O}(h^4)\,.
  \end{split}
\end{equation}

Similar identies can be stated for the finite volume discretizations
\begin{equation}
\begin{aligned}
\left(\Delta_\parallel^\VFT \bar \Phi^h\right)_{i,j} = \left(\Delta_\parallel^\ASYM \bar \Phi^h\right)_{i,j} \,,
\end{aligned}
\end{equation}

\begin{equation}
\begin{split}
\left(\Delta^\VFT \bar \Phi^h\right)_{i,j} &= \Delta \phi(x_i,y_j) \\
& \qquad +\frac{h^2}{12}\left( \frac{\partial^4}{\partial x^4} \phi(x_i,y_j)+ \frac{\partial^4}{\partial y^4} \phi(x_i,y_j)\right)+ \mathcal{O}(h^4)\,.
\end{split}
\end{equation}
\end{subequations}

\begin{subequations}
\begin{equation}
\begin{aligned}
&\left(\Delta_\parallel^\VFF \bar \Phi^h\right)_{i,j} = \Delta_\parallel \phi(x_i,y_j) \\
&\hspace*{5em}+\alpha^2 \frac{h^4}{1920}\left(\frac{\partial^6}{\partial x^2 \partial y^4} \phi(x_i,y_j) - \frac{64}{3}\frac{\partial^6}{\partial x^6} \phi(x_i,y_j) \right)\\
	&\hspace*{5em}- \alpha\beta\frac{h^4}{15}\left(\frac{\partial^6}{\partial x^5 \partial y} \phi(x_i,y_j) + \frac{\partial^6}{\partial x \partial y^5} \phi(x_i,y_j)\right) \\
	&\hspace*{5em} +\beta^2 \frac{h^4}{1920}\left(\frac{\partial^6}{\partial x^4 \partial y^2} \phi(x_i,y_j) - \frac{64}{3}\frac{\partial^6}{\partial y^6} \phi(x_i,y_j) \right) + \mathcal{O}(h^6)\,.
\end{aligned}
\end{equation}
\begin{equation}
\begin{aligned}
\left(\Delta^\VFF \bar \Phi^h\right)_{i,j} &= \Delta\phi(x_i,y_j) +\frac{h^4}{1920}\Bigg(\frac{\partial^6}{\partial x^2 \partial y^4} \phi(x_i,y_j) - \frac{64}{3}\frac{\partial^6}{\partial x^6} \phi(x_i,y_j) \\	& \qquad \quad + \frac{\partial^6}{\partial x^4 \partial y^2} \phi(x_i,y_j) - \frac{64}{3}\frac{\partial^6}{\partial y^6} \phi(x_i,y_j) \Bigg) + \mathcal{O}(h^6)\,.
\end{aligned}
\end{equation}
\end{subequations}

\end{document}